\documentclass[12pt,a4paper]{article}
\usepackage{amsmath,amscd,amssymb,amsfonts,amsthm}

\numberwithin{equation}{section}

\title{On complexes related with calculus of variations
\thanks{The research was supported in part by EPSRC, under grant GR/N00821/01.}}

\newcommand{\myaddress}{\begin{center}
\small Department of Mathematics\\ \small University of Manchester Institute of Science
and Technology (UMIST)\\ \small PO Box 88, Manchester M60 1QD, England\\ \small {\tt
theodore.voronov@umist.ac.uk, khudian@umist.ac.uk}
\end{center}}

\author{Hovhannes Khudaverdian, Theodore Voronov}

\date{\myaddress}
\renewcommand{\leq}{\leqslant}
\renewcommand{\geq}{\geqslant}

\newcommand{\dd}{{\,\,}{\bar{\smash{\!\!d}}}}

\newtheorem{prop}{Proposition}[section]

\theoremstyle{definition}
\newtheorem{de}{Definition}[section]

\newtheorem{rem}{Remark}[section]

\DeclareMathOperator{\Ber}{Ber}

\DeclareMathOperator{\sign}{sign}

\DeclareMathOperator{\Vect}{Vect} 
\newcommand{\lie}[1]{{\cal L}_{{#1}}}
\renewcommand{\P}[1]{{\cal P}_{{#1}}}
\newcommand{\volx}{{\cal D}{x}}
\newcommand{\volt}{{\cal D}{t}}
\newcommand{\der}[2]{{\frac{\partial {#1}}{\partial {#2}}}}
\newcommand{\lder}[2]{{\partial {#1}/\partial {#2}}}
\newcommand{\dder}[3]{{\frac{\partial^2 {#1}}{\partial {#2}\partial {#3}}}}
\newcommand{\Om}[2]{{\boldsymbol\O^{#1|#2}}}
\newcommand{\R}[1]{{\mathbb R}^{#1}}

\newcommand{\p}{\partial}

\newcommand{\fun}{C^{\infty}}

\def\a{\alpha}
\def\e{\varepsilon}
\def\s{\sigma}
\def\f{{\varphi}}
\def\O{\Omega}
\def\G {\Gamma}
\def\o{\omega}
\newcommand{\g}{{\gamma}}
\newcommand{\F}{{\Phi}}
\renewcommand{\S}{{\Sigma}}
\renewcommand{\r}{{\rho}}
\newcommand{\x}{{\xi}}
\def\d{\delta}

\renewcommand{\L}{{\bf L}}

\newcommand{\itt}{{\tilde\imath}}
\newcommand{\ot}{{\tilde \omega}}
\newcommand{\at}{{\tilde a}}
\newcommand{\bt}{{\tilde b}}
\newcommand{\mut}{{\tilde \mu}}
\newcommand{\Xt}{{\tilde X}}
\newcommand{\Yt}{{\tilde Y}}
\newcommand{\Kt}{{\tilde K}}

\newcommand{\fsm}{{\f_{\s}^{\mu}}}
\renewcommand{\fam}{{\f_{a}^{\mu}}}
\newcommand{\fasm}{{\f_{a\s}^{\mu}}}
\newcommand{\Gsm}{{\G_{\s}^{\mu}}}
\newcommand{\Bsm}{{B^{\s}_{\mu}}}
\newcommand{\osm}{{\o^{\s}_{\mu}}}
\newcommand{\Ysm}{{Y_{\s}^{\mu}}}
\newcommand{\Gasm}{{\G_{a\s}^{\mu}}}
\newcommand{\Yasm}{{Y_{a\s}^{\mu}}}
\newcommand{\fmu}{{{\cal F}_{\mu}}}
\newcommand{\fm}{{\f^{\mu}}}
\newcommand{\fa}{{{\cal F}_{a}}}
\newcommand{\xia}{{x_{i}^{a}}}

\newcommand{\xija}{{x_{ij}^{a}}}
\newcommand{\xja}{{x_{j}^{a}}}

\begin{document}
\maketitle

\begin{abstract}
We consider the variational complex on infinite jet space and the complex of variational
derivatives for Lagrangians of multidimensional paths and study relations between them.
The discussion of the variational (bi)complex is set up in terms of a flat connection in
the jet bundle. We extend it to supercase using a particular new class of forms. We
establish relation of the complex of variational derivatives and the variational complex.
Certain calculus of Lagrangians of multidimensional paths is developed.  It is shown how
covariant Lagrangians of higher order can be used to represent characteristic classes.
\end{abstract}

\tableofcontents

\section*{Introduction}

It is possible to include the ``Euler-Lagrange operator'' (the left-hand side of the
Lagrange equations) into a complex on the space of infinite jets, the so-called
variational complex. In this approach, Lagrangians (multiplied by volume forms) appear in
a particular term of this complex, the other terms consisting of certain forms or classes
of forms on the jet space. The explanation of this complex is in a spectral sequence due
to A.M.~Vinogradov. On the other hand, there is another complex, which we call the
complex of variational derivatives, (suggested by one of the authors), and in which the
variational derivative is the main ingredient of the differential for all terms. Each
term consists of Lagrangians of multidimensional paths, and the differential increases
the dimension of paths by one. The purpose of this paper is to initiate the study of the
relation of these two complexes connected with variational problems.

The paper consists of three sections. In Sections 1 and 2 we mainly review the material
concerning the variational complex. Obviously, it is well known to experts. However, we
make a point in systematically using the framework of the canonical flat connection in
the infinite jet bundle and following analogies with the more familiar
differential-geometric constructions. Interpreting the Cartan distribution (the contact
distribution) as a connection works  only for infinite jets. Hence, we make them, not
finite jets, our primary objects.  Hopefully, this exposition can help to straighten out
and clarify some points.

In Section 1, we define the Cartan connection $\G$ in the jet bundle. It gives the
exterior covariant differential on horizontal forms, which extends to the ``horizontal''
differential on arbitrary forms. The ``vertical'' differential is obtained as the Lie
derivative along $\G$ considered as a form-valued vector field. This gives a canonical
bicomplex on the jet space. In Section 2, Vinogradov's spectral sequence and variational
complex are described and the relation with the classical variational problem is
explained.

In Section 3, we introduce the complex of variational derivatives. Then the variational
complex is revisited and generalized for supermanifolds. Since the straightforward
generalization of the usual bicomplex is not satisfactory, we introduce a new class of
forms that are material for this case. These forms are particular hybrids of integral and
differential forms. Then we establish relation between the complex of variational
derivatives and the variational complex using a sequence of bicomplexes corresponding to
increasing dimension of the manifolds of parameters.  Then in the framework of the
complex of variational derivatives we study covariant Lagrangians of first and higher
order, which are natural integration objects over surfaces. We begin to develop a
``calculus of covariant Lagrangians'', in particular, we consider the composition of
Lagrangians and relate it with some densities for characteristic classes. This study is
in progress and we hope to elaborate it elsewhere. Potentially it can be linked with
symplectic reduction and various invariants for supermanifolds.

Main sources on jet geometry and variational complex are the books by Vinogradov and
others~\cite{vinog:kollektiv} and Olver~\cite{olver}. They contain plenty of reference
for other works. Experts might notice that though our exposition in the preliminary
sections owe much to these basic sources, our approach in many points differs from
both~\cite{vinog:kollektiv},~\cite{olver}. In our work we use supermathematics. Not only
we try to generalize to supercase, but we substantially rely on supermethods, which
simplify constructions important for the geometry of jets. For supermanifolds theory we
refer to~\cite{berezin:antieng}, two chapters from~\cite{manin:gaugeeng},
to~\cite{leites:book} and to~\cite{tv:git}. Concerning integration theory we particularly
refer to~\cite{tv:git} and~\cite{tv:dual}.

\smallskip
\textbf{Acknowledgements.} One of us (H.K.) is grateful to A.~Verbovetsky for many
inspiring discussions of jet geometry in 1998-99.

\section{Preliminaries: Cartan connection on infinite jet space}\label{secprelim}
\subsection{Infinite jet space}

Consider a fibre bundle $\pi=\pi(E,M,F)$ over an $r$-dimensional manifold $M$ with the
fibre $F$.

Consider the fibre bundle $J^k(\pi)\rightarrow M$ of $k$-jets of local sections of $E$.
If $x^a$ $=(x^1,\dots,$ $x^m)$ and ${\varphi}^\mu=({\varphi}^1,\dots,{\varphi}^n)$ are
some local coordinates on $M$ and $F$ respectively, then $(x^a,{\varphi}^\mu_\s)$
$=(x^a,{\varphi}^\mu,{\varphi}^\mu_{a_1},{\varphi}^\mu_{a_1 a_2},$
$\dots,{\varphi}^\mu_{a_1\dots a_k})$ are natural local coordinates on $J^k(\pi)$. Here
$\s=a_1\dots a_p$ is a multi-index, $p=|\s|\leq k$. Notice that the natural bundles
$J^k(\pi)\rightarrow E$, $J^k(\pi)\rightarrow J^l(\pi)$ ($k>l$) are not vector bundles,
though their respective fibres are Euclidean spaces.

Consider the space  $J^{\infty}(\pi)$ of infinite jets, i.e., the inverse limit of the
manifolds $J^k(\pi)$. We denote by $C^\infty(J^k(\pi))$ the space of smooth functions on
$J^k(\pi)$ $(k=0,1,2,\dots,\infty)$. Every function $f$ on $J^\infty(\pi)$, by the
definition, has finite order, i.e., depends on a finite number of variables:
$f=f(x,({\varphi}^\mu,\dots,{\varphi}^\mu_{\a_1\dots\a_p}))$.  In other words, $f\in
C^\infty(J^p(\pi))$ for some finite $p$.

To every local section $s(x)$ of the fibre bundle $E$ corresponds its $k$-jet, $(j_k
s)(x)$ $=(x^a,\f^\mu(x),\f^\mu_a(x),\dots,\f^\mu_{a_1\dots a_k}(x))$, where
$\f^\mu_{a_1\dots a_p}(x)=$ ${\p^p \f^\mu(x)\over\p x^{a_1}\dots \p x^{a_p}}$, which is a
section of $J^k(\pi)$, and the infinite jet $(j s)(x)=$ $(x^a, \f^\mu(x),$ $\f^\mu_a(x),
\dots,$ $\f^\mu_{a_1\dots a_p}(x),\dots)$, which is a section of $J^{\infty}(\pi)$.

In the sequel we shall denote $J:=J^{\infty}(\pi)$. Infinite jets will be simply called
``jets''. We shall use the notation $[\f]$ for the whole collection $\fsm$.

For a function $f=f(x,[\f])$ on $J$ and a local section $s(x)$ of  the bundle
$E\rightarrow M$ the value $f\vert_s$ is defined as
\begin{equation}\label{functionvalue}
  f\big\vert_{s}=f\circ js =f\left(x,[\f(x)]\right).
\end{equation}
This is a function on $M$.

Vectors on the space of infinite jets $J$ are derivations of the algebra $C^\infty(J)$.
They have the form
\begin{equation}\label{vector}
  X=X^a{\p\over \p x^a}+
    X^\mu_\s{\p\over \p {\varphi}^\mu_\s},
\end{equation}
where  the number of non-zero coefficients can be infinite. (Summation in~(\ref{vector})
is over multi-indices of all orders.) A vector on $J$ is \textit{vertical} if  its
projection on $M$ vanishes, i.e., $X^a$ in~\eqref{vector} equals to zero.

Consider now the algebra $\O=\O(J)$  of differential forms on the space of infinite jets.
Differential forms that vanish on vertical vectors  are called \textit{horizontal} forms.
They have the appearance $\omega= \sum \omega_{a_1\dots a_p}(x,{\varphi}^\mu_\s)\,
dx^{a_1}\dots dx^{a_1}$. (We don't write the wedge sign. Instead, we use a convention of
supermathematics that the differential $dx^a$ is odd for an even variable $x^a$.)

\subsection{The Cartan connection}

In the fibre bundle $J\to M$ there is a natural connection specified by the vector-valued
one-form
\begin{equation}\label{cartanconnection}
  \G=  \sum_{\mu,\s}
           \G^\mu_\s
           {\p\over\p {\varphi}^\mu_\s},
\end{equation}
which takes values in vertical vectors. (Summation in~(\ref{cartanconnection}) goes over
all multi-indices.) Here the coefficients
\begin{equation}\label{cartanconnectioncom}
  \G_\s^\mu= d{\varphi}^\mu_\s- dx^b\,{\varphi}^\mu_{b\s}
\end{equation}
are the so called Cartan one-forms. We call $\G$ the \textit{Cartan connection}. Here $b
\s$ is the multi-index $b a_1\dots a_k$, if $\s$ is the multi-index $a_1\dots a_k $.

The distribution of the horizontal planes w.r.t. the connection $\G$ is called the
\textit{Cartan distribution}, its dimension is $m={\rm dim}M$.  The ideal in the algebra
$\O(J)$ corresponding to the Cartan distribution (consisting of all forms vanishing on
the distribution) will be denoted $C\O\subset \O$. It is generated by the Cartan forms
\eqref{cartanconnectioncom}. (Notice that the Cartan distribution is defined on
$J^k(\pi)$ for finite $k$ also, but  there it does not correspond to a connection.)

\textit{The Cartan connection in the bundle $J\to M$ is flat.} That means that the ideal
$C\O$ is a differential ideal: $d(C\O)\subset C\O$. The flatness of the Cartan connection
(the integrability of the Cartan distribution) is an essential feature of the infinite
jet space. For finite $k$ the Cartan distribution on $J^k(\pi)$ is not integrable.

Every tangent vector on $J$ can be uniquely decomposed into a vertical and a horizontal
vector: $X=X_{\rm vert}+X_{\rm hor}$, where $X_{\rm vert}$ $:=\langle X,\G \rangle$, the
value of the one-form $\G$ on $X$, and $X_{\rm hor}:=X-X_{\rm vert}$. $X_{\rm hor}$
belongs to the Cartan distribution, $\langle X_{\rm hor},\G\rangle=0$:

\begin{equation}\label{vectordecomposing}
  X=
  X^a{\p\over \p x^a}+X^\mu_\s{\p\over\p {\varphi}^\mu_\s}=
     \underbrace {X^a
     \left(
      {\p\over \p x^a}+
             {\varphi}^\mu_{a\s}{\p\over\p {\varphi}^\mu_\s}
             \right)}_{\rm horizontal}+
      \underbrace
      {\left(X^\mu_\s-X^a {\varphi}^\mu_{a\s}\right)
      {\p\over\p {\varphi}^\mu_\s}}_
      {\rm vertical}
\end{equation}

Let us emphasize that every connection in an arbitrary fibre bundle defines two
operations: a \textit{covariant derivative of sections} and a \textit{covariant
derivative of functions on the total space}. If $X$ is a vector on the base at a point
$x_0$, then for a section $s$, its covariant derivative along $X$, $D_Xs$, is a vertical
tangent vector to the total space at the point $s(x_0)$, the vertical component of the
vector $s_*X$; for a function $f$ on the total space, its covariant derivative along $X$,
$D_Xf$, at a point $y_0$ of the total space, is the usual derivative of $f$ along the
horizontal lift of the vector $X$. (For the familiar case of vector bundles, these
notions correspond to the usual covariant derivative of sections, and to that of sections
of the dual bundle, which can be treated as linear functions on fibres.)

In particular, for the bundle $J\to M$ we obtain for the covariant derivative of a
section  $j(x)$ $=(x^a,{\varphi}^\mu(x),{\varphi}^\mu_a(x), {\varphi}^\mu_{ab}(x),\dots)$
\begin{equation}\label{sectionsderiv}
  D_Xj=\langle j_*(X),\G\rangle=
                     X^a
           \left(
              \der{{\varphi}^\mu_\s}{x^a} -{\varphi}^\mu_{a\s}
           \right)
           {\p\over\p {\varphi}^\mu_\s}.
\end{equation}
For a function $f$ on $J$ we obtain
\begin{equation}\label{functionsderiv}
  D_Xf=X^a D_af=
               X^a
             \left(
         {\p f\over \p x^a}+
      {\varphi}^\mu_{a\s}{\p f\over\p {\varphi}^\mu_\s}
             \right).
\end{equation}
The covariant derivatives of functions  and sections  are related by the following analog
of the Leibniz formula:
\begin{equation}\label{leibnitz}
  X(j^*f)=(D_X j)f + j^* (D_Xf),
\end{equation}
where $f\in\fun(J)$, $j: M\to J$ is a section, and $X=X^b{\p/\p x^b}$ is a vector field
on $M$. (Notice that $D_X j\in \fun (M, j^*TJ)$ is a vector field along the map $j$.)

The jet $js(x)$ of an arbitrary section $s(x)$ of fibre bundle $E\to M$ gives an integral
surface of the Cartan distribution. It is clear from~\eqref{sectionsderiv} that  $D_X
js=0$ for all vector fields $X$. Hence, from~\eqref{sectionsderiv}\,--\,\eqref{leibnitz}
it follows that
\begin{equation}\label{reasonofunderstanding}
  X (f\big\vert_s)=D_Xf\big\vert_s
\end{equation}
where $f\big\vert_s$ is the value of $f$ on $s$, as defined above
in~\eqref{functionvalue}.

The covariant derivative  allows to define the \textit{covariant differential} of a
function on the jet space as a horizontal $1$-form:
\begin{equation}\label{functionsdiff}
  Df=dx^aD_af =
            dx^a\left(
          {\p f\over\p x^a}+
            {\varphi}^\mu_{a\s }
            {\p f\over\p {\varphi}^\mu_\s}
            \right).
\end{equation}
The operation $D:\fun(J)\to \O^1_{\text{hor}}(J)$ extends to the \textit{exterior
covariant differential} on horizontal forms of arbitrary degree:
\begin{equation}
  D\left(\o_{a_1\dots a_k}dx^{a_1}\dots dx^{a_k}\right)=
\left(D\o_{a_1\dots a_k}\right) dx^{a_1}\dots dx^{a_k}).
\end{equation}
In a coordinate-free notation,
\begin{equation}\label{formsdiff}
  D\o=d\o-\G\o,
\end{equation}
where we denote
\begin{equation}\label{firstoperation}
  \G\o=  \G^\mu_\s{\p\o_{a_1\dots a_k}\over\p {\varphi}^\mu_\s}\,
       dx^{a_1}\dots dx^{a_k}
\end{equation}
(this makes sense only for horizontal forms). We arrive at the sequence
\begin{equation}\label{horizontalcomplex}
 0\to \O_{\rm hor}^0{\buildrel D\over\longrightarrow}\,
     \O_{\rm hor}^1{\buildrel D\over\longrightarrow}\,
                   \dots
        {\buildrel D\over\longrightarrow}\,
        \O_{\rm hor}^{m-1}{\buildrel D\over\longrightarrow}\,
                  \O_{\rm hor}^m\to 0,
\end{equation}
where $\O_{\rm hor}^k$ is the space of horizontal $k$-forms. Due to flatness of the
Cartan connection,  $D^2=0$.  Hence,~\eqref{horizontalcomplex} is a complex. This complex
is the first part of the variational complex defined in the next section.

\subsection{The bicomplex $\O^{**}(J)$} \label{subsecbicomplex}

The operation~\eqref{firstoperation} and the differential~\eqref{formsdiff} introduced
above can be extended from horizontal forms to the whole algebra $\O(J)$. This can be
done using  operations with form-valued vector fields. Let us briefly recall the
necessary notions.

It is very convenient to use ``super'' language. A differential form on a manifold $M$
can be considered as a function of even (commuting) variables $x^a$ and odd
(anticommuting) variables $dx^a$: $\o=\o(x,dx)$. We denote parity of all objects by
tilde: $\tilde x^a=0$, $\tilde dx^a=1$.  The variables $x^a, dx^a$ are coordinates on the
supermanifold ${\Pi}TM$ associated with tangent bundle $TM$ of $M$. In these terms the
exterior differential $d$ on $\O(M)$ is nothing but the vector field $dx^i\,{\p/ \p
x^i}\in\Vect({\Pi}TM)$.

A form-valued vector field $X\in\Vect(M, \O(M))$ has the appearance
\begin{equation}
  X=X^i(x,dx)\der{}{x^i}.
\end{equation}
Warning: it is \textit{not} a vector field on $\Pi TM$. The \textit{Lie derivative ${\cal
L}_X$ along $X$} is the vector field on ${\Pi}TM$ defined by the conditions: ${\cal L}_Xf
=Xf$ for an arbitrary function on $M$ and ${\cal L}_X$ commutes with $d$:
\begin{equation}\label{supercont}
  [{\cal L}_X,d]={\cal L}_X\circ d-(-1)^{\Xt}
                      d\circ {\cal L}_X=0.
\end{equation}
Hence
\begin{equation}\label{explicitsupercont}
  {\cal L}_X=X^i(x,dx){\p\over\p x^i}+(-1)^{\Xt}
           dX^i(x,dx) {\p\over\p dx^i}.
\end{equation}
All vector fields on $\Pi TM$ commuting with $d$ have the appearance ${\cal L}_X$ for
some form-valued vector field $X$ on $M$. The \textit{Nijenhuis bracket} $[X,Y]$  of
form-valued vector fields $X$ and $Y$  is defined by the formula:
\begin{equation}\label{nijenhuis}
  \left[{\cal L}_X, {\cal L}_Y \right]={\cal L}_{[X,Y]}.
\end{equation}
It extends the usual commutator. Explicitly,
\begin{equation}\label{nijenhuis2}
  [X,Y]=\left({\cal L}_X Y^i -(-1)^{\Xt\Yt}{\cal L}_Y X^i\right)\der{}{x^i}.
\end{equation}

Now we can apply these constructions to our situation. The connection $1$-form $\G$ can
be viewed as a form-valued vector field on $J$:
\begin{equation}
  \G=\G_\s^\mu\der{}{\f^\mu_\s}=\left(d{\varphi}^\mu_\s- dx^a{\varphi}^\mu_{a\s}\right) \der{}{\f^\mu_\s}.
\end{equation}
The flatness of $\G$ is equivalent to
\begin{equation}\label{connectionflat}
  [\G,\G]=0
\end{equation}
(Nijenhuis bracket). It follows that
\begin{equation}
  [{\cal L}_{\G},{\cal L}_{\G}]=2\, \lie{\G}^2=0
\end{equation}
for the Lie derivative along $\G$.  From the explicit formula above, we get for  an
arbitrary form $\o\in\O(J)$:
\begin{equation}
  {\cal L}_{\G}\o= \left(d\f_\s^\mu - dx^a\f_{a\s}^\mu\right)\der{\o}{\f_\s^\mu}
  -dx^ad\f_{a\s}^\mu \der{\o}{d\f_\s^\mu}.
\end{equation}
Hence, ${\cal L}_{\G}x^a=0$, ${\cal L}_{\G}\f_\s^\mu=\G_\s^\mu$, ${\cal L}_{\G}dx^a=0$,
${\cal L}_{\G}d\f_\s^\mu=-dx^a d\f_{a\s}^\mu$. Notice that for horizontal forms ${\cal
L}_{\G}\o={\G}\o$ (with $\G\o$ being defined in~\eqref{firstoperation}). It also follows
that
\begin{equation}
  {\cal L}_{\G}\G_\s^\mu=0,
\end{equation}
because $[{\G},{\G}]=2({\cal L}_{\G}\G_\s^\mu)\,\lder{}{\f^\mu_\s}$. It turns out to be
convenient to express forms in terms of the variables $x^a$, $\f_\s^\mu$, $dx^a$,
$\G_\s^\mu$ instead of $x^a$, $\f_\s^\mu$, $dx^a$, $d\f_\s^\mu$. Written in these
variables, ${\cal L}_{\G}$ is simply
\begin{equation}
  {\cal L}_{\G}= \G_\s^\mu \der{}{\f_\s^\mu}
\end{equation}
for all forms.

Define on arbitrary forms two operations:
\begin{align}%
  \d \o& := {\cal L}_{\G}\o,  \label{vertderivaction} \\
  D\o & := d\o- {\cal L}_{\G}\o, \label{horizontalderivaction}
\end{align}
and call $\d$ and $D$ the \textit{vertical} and \textit{horizontal differentials},
respectively. Hence, $d=D+\d$. We know that $\d^2={\cal L}_{\G}^2=0$ and that $[d,
\d]=[d,{\cal L}_{\G}]=0$. Hence, $D^2=0$, and $[D,\d]=D\circ \d+\d\circ D=0$, so we have
a bicomplex. Clearly, the horizontal differential $D$ extends the covariant differential
of horizontal forms defined above.

We can introduce an invariant \textit{bigrading} in the algebra $\O(J)$ by the degrees in
the variables  $dx^a$ and $\G_{\s}^\mu$. By definition,
\begin{equation}\label{bigrading}
  \O^{p,q}(J):=\Bigl\{\o\in\O(J)\,\bigl\vert\, \#\Gsm(\o)=p \text{  \,and\,  } \#dx^a(\o)=q
  \bigr.\Bigr\},
\end{equation}
where $\#dx^a$, etc., mean the degree in the respective variables. The total degree is
$p+q$. Since one can easily obtain $D(x^a)=dx^a$, $D(\f_{\s}^{\mu})=dx^a \f_{a\s}^\mu$,
$D(dx^a)=0$, $D(\G_{\s}^\mu)=dx^a\G_{a\s}^\mu$ (from the formulas for $d$ and for
$\d={\cal L}_{\G}$ above), we have for a form $\o=\o(x,[\f],dx, \G)$:
\begin{align}
  D\o & =dx^a\left(\der{\o}{x^a}+\f_{a\s}^\mu \der{\o}{\f_{\s}^\mu}+ \G_{a\s}^\mu
  \der{\o}{\G_{\s}^\mu}\right), \label{horderivaction2} \\
  \d\o & = \G_\s^\mu \der{\o}{\f_\s^\mu}.\label{vertderivaction2}
\end{align}
Hence, $D: \O^{p,q}(J)\to\O^{p,q+1}(J)$, $\d: \O^{p,q}(J)\to\O^{p+1,q}(J)$.

(There is a remote analogy with complex manifolds, where the integrability is provided by
the condition $[J,J]=0$  and  there is a bigrading of forms by $dz$ and $d\bar z$. Of
course, in our case there is no such symmetry that exists between the holomorphic  and
antiholomorphic differentials.)

\subsection{Evolutionary vector fields}\label{subsecevol}

Consider the action of vector fields  on the connection $1$-form $\G$. Notice first that
geometrically $\G$ has the meaning of a projector onto the vertical subspace in a tangent
space to $J$, the kernel of this projection being the horizontal subspace (thus defined).
It is easy to see that for an arbitrary projector $P$, its preservation by some flow is
equivalent to the preservation of two distributions, the image of $P$ and the kernel of
$P$ (the image of $1-P$). The preservation of the image as such is equivalent to
$(1-P)\circ \lie{X}P=0$, while the preservation of the kernel is equivalent to $P\circ
\lie{X}P=0$. So in our case, vector fields that preserve $\G$, preserve both the Cartan
distribution and the fibre structure in $J\to M$. Preservation of the Cartan distribution
as such is equivalent to $\langle \lie{X}\G,\G\rangle=0$.

In the same way as we did with $1$-forms above, decompose the tangent space at some point
$j\in J$ into the horizontal subspace with a natural basis
$D_a=\lder{}{x^a}+\f_{a\s}^{\mu}\lder{}{\f_{\s}^{\mu}}$ and the vertical subspace with a
natural basis $\lder{}{\f_{\s}^{\mu}}$. We have: $\langle D_a, dx^b\rangle=\d_a^b$,
$\langle D_a, \G_{\s}^{\mu}\rangle=0$, $\langle \lder{}{\f_{\s}^{\mu}}, dx^b\rangle=0$,
$\langle \lder{}{\f_{\s}^{\mu}}, \G_{\tau}^{\nu}\rangle=\d_{\tau}^{\s}\d_{\mu}^{\nu}$.
Hence, this is the dual basis for the basis of $1$-forms $dx^a$, $\G_{\s}^{\mu}$.

It makes sense to consider horizontal and vertical vector fields separately.

For an arbitrary vector field $X\in\Vect (J)$,  the  Lie derivative of the Cartan form
$\G$ along $X$ is nothing but the Nijenhuis bracket: $\lie{X}\G=[X,\G]$. It follows that
at no extra cost we can include in our consideration vector fields that take values in
forms. (This might be useful, e.g., for studying deformations.) Assuming this, we have
$[X,\G]=-(-1)^{\Xt}[\G,X]$.

Consider first a (form-valued) horizontal vector field:
\begin{equation}
  X=X^a(x,[\f], dx, \G)\,D_a=X^a(x,[\f], dx,
  \G)\left(\der{}{x^a}+\fasm\der{}{\fsm}\right).
\end{equation}
To it corresponds the following infinitesimal transformation:
\begin{align*}
  x^a & \mapsto x^a+\e X^a\\
  \fsm & \mapsto \fsm+\e X^a\fasm \\
  dx^a & \mapsto dx^a+\e(-1)^{\Xt} dX^a\\
  d\fsm & \mapsto d\fsm+\e \left((-1)^{\Xt} dX^a\fasm +(-1)^{\Xt} X^a d\fasm\right).\\
\intertext{Hence, by a direct computation we obtain for $\Gsm=d\fsm-dx^a\fasm$:}
  \Gsm &\mapsto \Gsm+ \e X^a\Gasm.
\end{align*}
So, $\lie{X}\Gsm=X^a\Gasm$. Using the formula~\eqref{nijenhuis2} for the Nijenhuis
bracket, we obtain from here that
\begin{equation}
  \lie{X}\G=-(-1)^{\Xt}\d X^a\,D_a,
\end{equation}
where $\d$ is the vertical differential~\eqref{vertderivaction2}, or:
\begin{equation}
  [\G,X]= \d X^a D_a.
\end{equation}
Hence, $\langle\lie{X}\G,\G\rangle=0$ is trivially satisfied (as it should be),  while
$\lie{X}\G=0$ is equivalent to $\d X^a=0$ componentwise. In particular, for a usual
vector field, the condition $\d X^a=0$ simply means that $X^a=X^a(x)$ (is independent of
$\fsm$). We have proved the following statement (where we slightly changed our notation):
\begin{prop}\label{horpreserv}
Horizontal vector fields on $J$ preserve the Cartan distribution automatically; they
preserve also the bundle structure $J\to M$ if they are covariant derivatives
\begin{equation}
  D_X=X^a(x)\,D_a=X^a(x)\left(\der{}{x^a}+\fasm\der{}{\fsm}\right),
\end{equation}
where $X=X^a(x)\,\lder{}{x^a}\in\Vect(M)$. Clearly,
\begin{equation}
  [D_X,D_Y]=D_{[X,Y]},
\end{equation}
because the connection is flat.
\end{prop}

Notice that for general form-valued horizontal vector fields we obtain ``vertical''
complexes $\O^{*,q}(J,\Vect_{\text{hor}})\stackrel{\d}{\longrightarrow}
\O^{*+1,q}(J,\Vect_{\text{hor}})$ with the differential $\d$ applied componentwise.
Obviously, there are similar complexes with coefficients in arbitrary vector bundles
coming from $M$, i.e., with $\d$-flat transition functions  (independent of $\fsm$).

Consider now a vertical vector field
\begin{equation}
  Y=\Ysm(x,[\f],dx,\G)\der{}{\fsm}
\end{equation}
(taking values in forms). 
Again, we can directly find that $\lie{Y}\Gsm=(-1)^{\Yt}(d\Ysm-dx^a\Yasm)$. Hence, it
follows that
\begin{equation}\label{cartanformlieder}
  \lie{Y}\G=[Y,\G]=(-1)^{\Yt}dx^a\left(D_a Y^\mu_\s-Y^\mu_{a\s}\right)\der{}{\fsm},
\end{equation}
or:
\begin{equation}\label{hordiffonfields}
  [\G,Y]=-(-1)^{\Yt}[Y,\G]=\left(dx^a\Yasm-D \Ysm\right)\der{}{\fsm}.
\end{equation}
The formula~\eqref{hordiffonfields} can be interpreted as defining the differential in a
complex of form-valued vertical vector fields. For a usual vertical vector field, we
deduce that the condition $\lie{Y}\G=\langle\lie{Y}\G,\G\rangle=0$ gives rise to the
inductive formula
\begin{equation}
  \Yasm=D_a\Ysm,
\end{equation}
hence all coefficients $\Ysm$ are expressed from the coefficient
$Y^{\mu}=Y^{\mu}(x,[\f])$ as
\begin{equation}\label{evolvectorfields}
  \Ysm=D_{\s}^{|\s|}Y^{\mu}=D_{a_1}\dots D_{a_k}Y^{\mu},
\end{equation}
where $k=|\s|$ is the order of the multi-index $\s=a_1\dots a_k$. In other words, a
vertical vector field  preserving the Cartan connection is uniquely defined by the
initial term $Y^{\mu}\lder{}{\f^{\mu}}$, which is the restriction of the field to the
subalgebra of functions of $x^a,\f^{\mu}$.

Changing the notation, let $Y$ be a ``generalized'' vertical vector field on the bundle
$E$ with coefficients depending on jets:
\begin{equation}    \label{evolutionary}
  Y=Y^{\mu}(x,[\f])\,\der{}{\f^{\mu}}.
\end{equation}
More formally, fields like~\eqref{evolutionary} are sections of the pull-back to $J$ of
the vertical subbundle of $TE$. Denote by $\P{Y}$ the vertical vector field on $J$
uniquely defined by the conditions: $\P{Y}f=Yf$ for all functions on $E$ and
$\lie{\P{Y}}\G=[\P{Y},\G]=0$. Explicitly,
\begin{equation}\label{evolutionary2}
  \P{Y}=\sum D_{\s}^{|\s|}Y^{\mu}\der{}{\fsm}.
\end{equation}
``Generalized'' vector fields of the form~\eqref{evolutionary} on the bundle $E$ (to
which correspond genuine vector fields $\P{Y}$ on $J$) are called \textit{evolutionary
vector fields}. They model the variations of the independent argument (variations of
sections) of the classical calculus of variations. See below in Section~\ref{secvin}. We
have obtained the following statement:

\begin{prop}\label{vertpreserv}
All vertical vector fields on $J$ preserving the Cartan distribution have the form
$\P{Y}$ for some evolutionary vector field $Y$. Since they form a subalgebra in
$\Vect(J)$, the formula
\begin{equation}\label{jacobi}
  [\P{Y},\P{Z}]=\P{[Y,Z]}
\end{equation}
defines a bracket $[Y,Z]$ on evolutionary vector fields, which extends the usual
commutator.
\end{prop}
(Notice an analogy with the construction of the Nijenhuis bracket.) The bracket of
evolutionary vector fields defined by~\eqref{jacobi} is known as the \textit{Jacobi
bracket}. Explicitly, for evolutionary vector fields $Y$ and $Z$,
\begin{equation}\label{jacobiexplicit}
  [Y,Z]=(\P{Y}Z^{\mu}-\P{Z}Y^{\mu})\der{}{\f^{\mu}}.
\end{equation}
Notice the following useful properties of evolutionary fields:
\begin{align}
  [D,\iota_{\P{Y}}]&=0, \label{Dwithev} \\
  [\d,\iota_{\P{Y}}]&=\lie{\P{Y}}, \label{deltawithev}
\end{align}
where $\iota$ stands for the interior multiplication. Indeed, the
equality~\eqref{Dwithev} follows by a direct calculation from~\eqref{horderivaction2} and
\eqref{evolutionary2}, and the equality~\eqref{deltawithev} follows then
from~\eqref{Dwithev} and the usual relation between $d$ and $\iota$. (Alternatively,
these formulas can be deduced using the extension of Cartan's identities to form-valued
vector field.) It also follows  that
\begin{equation}\label{lieevwithD}
  [\lie{\P{Y}},D]=0
\end{equation}
for all evolutionary fields $Y$.

Bringing the results of Propositions~\ref{horpreserv} and~\ref{vertpreserv} together, we
conclude that every vector field on $J$ preserving the Cartan connection $\G$ can be
uniquely decomposed into
\begin{equation}\label{genfieldprescart}
  D_X+\P{Y},
\end{equation}
where $X$ is a vector field on $M$ and $Y$ is an evolutionary vector field on $E$. If we
do not require that the bundle structure $J\to M$ is preserved, then $D_X$
in~\eqref{genfieldprescart} can be replaced by an arbitrary horizontal field. It is easy
to see that
\begin{equation}
  [D_X,\P{Y}]=0
\end{equation}
for all $X\in\Vect(M)$ and evolutionary $Y$, and that arbitrary horizontal fields form an
ideal in the Lie algebra of vector fields preserving the Cartan distribution.

\begin{rem}
Every vector field on the total space $E$,
\begin{equation}\label{fieldone}
  X=X^{a}(x,\f)\der{}{x^a}+X^{\mu}(x,\f)\der{}{\f^{\mu}},
\end{equation}
acts on sections $s:M\to E$ by the Lie derivative (infinitesimal variation):
\begin{equation}\label{liederofsectionaction}
  \bigl(\lie{X}s\bigr)(x)=\left(X^{a}(s(x))\,\der{\f^{\mu}}{x^a}(x)
  -X^{\mu}(s(x))\right)\der{}{\f^{\mu}},
\end{equation}
which is a vector field along the map $s:M\to E$. We can reinterpret this variation as
coming from an evolutionary vector field of a particular form:
\begin{equation}\label{liederofsection}
  \lie{X}:=\bigl(X^{a}(x,\f)\fam-X^{\mu}(x,\f)\bigr)\der{}{\f^{\mu}}.
\end{equation}
The subalgebra of such fields under the Jacobi bracket is isomorphic to the Lie algebra
$\Vect(E)$.
\end{rem}

\begin{rem}
To a vector field $X\in\Vect(E)$  corresponds its ``prolongation''
$X^{(\infty)}\in\Vect(J)$ defined by the conditions that $X^{(\infty)}f=Xf$ for all
functions on $E$ and that $X^{(\infty)}$ preserves the Cartan distribution. Decomposing
$X^{(\infty)}$ into the horizontal and vertical components, it is not difficult to see
that for the field~\eqref{fieldone} we get
\begin{equation}\label{xinfinity}
  X^{(\infty)}=X^{a}(x,\f)D_a-\P{Y},
\end{equation}
where $Y=\lie{X}$ is the evolutionary vector field~\eqref{liederofsection} corresponding
to $X$. (The minus sign before $\P{Y}$ in~\eqref{xinfinity} is explained by a ``duality''
of functions and sections, hence the opposite signs in Lie derivatives.) Notice, finally,
that the space of all fields on $J$ preserving the Cartan distribution is substantially
larger than the space of such prolongations.
\end{rem}

\section{Vinogradov's spectral sequence and variational complex}\label{secvin}
\subsection{Main construction}

In the complex $(\O, d)$ of differential forms on infinite jet space $J$ there is a
natural filtration by the powers of the Cartan ideal $C\O$:
\begin{equation}\label{vinfiltration}
   \ldots\subset C^k\O\subset C^{k-1}
            \subset\ldots\subset\O\subset C\O\subset C^0\O=\O,
\end{equation}
where $C^k\O$ denotes the $k$-th power of the ideal $C\O$ and $d$ is the usual de Rham
differential. $d(C^k\O)\subset C^k\O$, because the Cartan distribution is integrable. By
standard homological algebra, we come to the \textit{Vinogradov spectral
sequence}~\cite{vinog:field},
\begin{equation}\label{specsequence}
  (E_r^{**},d_r) \Longrightarrow H^*_{DR}(J)
\end{equation}
converging to the de Rham cohomology of the space $J$, with the zeroth term
\begin{equation}\label{vinzeroterm}
  E_0^{p,q}=C^p\O^{p+q}/C^{p+1}\O^{p+q}.
\end{equation}
In our case of a bundle $E\to M$ (see Remark~\ref{projective} below) the
filtration~\eqref{vinfiltration} is actually induced by the bigrading~\eqref{bigrading}:
\begin{equation}
  C^p\O^{p+q}=\O^{p,q}\oplus\O^{p+1,q-1}\oplus\O^{p+2,q-2}\oplus\ldots,
\end{equation}
and the Vinogradov spectral sequence is a spectral sequence of the bicomplex
$(\O^{**}(J), D, \d)$. We can identify its zeroth term~\eqref{vinzeroterm} with
$\O^{**}$:
\begin{equation}\label{vinzeroterm2}
  E_0^{p,q}=\O^{p,q}(J)
\end{equation}
and the differential $d_0$ is the horizontal differential $D:\O^{p,q}\to\O^{p,q+1}$. The
zeroth row of $E_0^{**}$ is the complex~\eqref{horizontalcomplex} of horizontal
differential forms on $J$.

Clearly, $\O^{p,q}=0$ for $q\geq m+1$ ($m$ is the dimension of the base $M$). Since the
fibres of the bundle $J\to E$ are contractible, we may say that the spectral sequence
converges to $H^*(E)$. The following essential fact holds:

\begin{prop}
If $p\geq 1$, then for all  $q<m$
\begin{equation}
  E_1^{p,q}=H^q(\O^{p,*},D)=0.\label{purealgebraic}
\end{equation}
Hence,   $E^{0,q}_1=E^{0,q}_\infty=H^q(E)$ for $q\leq m-1$, and
$E^{p,m}_2=E^{p,m}_\infty=H^{p+m}(E)$ for $p\geq 0$.
\end{prop}

This fact has  a purely algebraical origin [Vin]. We shall illustrate the main idea of
the proof of this Proposition when analyzing the content of the space $E^{1,m}_1$ (see
below).

Because of the degeneration property~\eqref{purealgebraic}, all the information about the
spectral sequence~\eqref{specsequence} is contained in the following
complex~\cite{vinog:kollektiv}:
\begin{equation}\label{varcomplex}
  E^{0,0}_0  {\buildrel d_0\over\longrightarrow}
   E^{0,1}_0  {\buildrel d_0\over\longrightarrow}\dots
              {\buildrel d_0\over\longrightarrow}
   E^{0,m-1}_0{\buildrel d_0\over\longrightarrow}
   E^{0,m}_0  {\buildrel d_1\circ p\over\longrightarrow}
   E_1^{1,m}{\buildrel d_1\over\longrightarrow}
   E_1^{2,m}{\buildrel d_1\over\longrightarrow}
   E_1^{3,m}{\buildrel d_1\over\longrightarrow}\dots
\end{equation}
where $p$ is the projection $E^{0,m}_0{\buildrel p\over\longrightarrow}E_1^{0,m}$.
This complex is called the \textit{variational
complex}.

Two halves of the complex~\eqref{varcomplex} can be described in terms of the bicomplex
$\O^{**}(J)$ as follows. The first half is the complex of horizontal
forms~\eqref{horizontalcomplex}: $E_0^{0,q}=\O^{0,q}$ with $d_0=D$; the second half
consists of classes of forms: $E_1^{p,m}=\O^{p,m}/D(\O^{p,m-1})$, and $d_1$ is the
vertical differential $\d$ acting on classes.

Hence, the cohomology of the complex~\eqref{varcomplex} in dimensions $0,\ldots, m-1$
coincides with $E_1^{0,k}$, $k=0,\ldots,m-1$. In dimensions $k\geq m$ the cohomology
coincides with $E_2^{0,k}$. It follows that the cohomology of the variational complex in
all dimensions is exactly $H^*_{DR}(J)=H^*_{DR}(E)$. In particular, the variational
complex is acyclic  after $(m+n)$-th term (where $n$ is the dimension of the fibre).

\begin{rem}\label{projective}
In a more general case of the so-called ``projective
jets''~\cite{vinog:kollektiv},\cite{olver} for an \mbox{$(m+n)$-dimensional} manifold
$E$, without a bundle structure, there is no bicomplex, but the
filtration~\eqref{vinfiltration} and the Vinogradov spectral sequence survive. The space
of projective jets is not contractible to $E$.
\end{rem}

\subsection{Relation with the classical variational problem}
\label{subsecvinclassicrelation}

The first part of the variational complex~\eqref{varcomplex} is the complex of horizontal
differential forms~\eqref{horizontalcomplex}. Horizontal forms of top degree ($m$-forms)
${\bf L}=L(x,[\f])\volx$ are \textit{Lagrangians}. We denote by $\volx$ the coordinate
volume form $dx^1\dots dx^m$. Later we shall also use the notation $\volx_a$ for the
$(m-1)$-form $(-1)^{a-1}dx^1\dots dx^{a-1} dx^{a+1}\dots dx^m$. If $s(x)=(x,\f(x))$ is a
section of the fibre bundle $E\to M$, then the value of the Lagrangian on this section
${\bf L}\big\vert_{s}=L\big\vert_{s}\volx$ (see~\eqref{functionvalue}) defines a
top-degree differential form ($m$-form) on the manifold $M$. Hence a Lagrangian
$\mathbf{L}$ defines an action functional  on sections  of $E$:
\begin{equation}\label{actionfunctional}
  S[\f]=\int_M {\bf L}=\int_M L\left(x,\f(x),\der{\f}{x}(x),\dots\right)\volx.
\end{equation}
Solution of the variational problem for this functional leads to the Euler-Lagrange
equations for the section $s(x)$:
\begin{equation}
  \fmu\big\vert_{s}=0,
\end{equation}
where the \textit{variational derivative} $\fmu=\fmu(\bf L)$ is defined by the following
expression
\begin{equation}\label{varderivative}
  \fmu({\bf L})=\der{L}{\fm}-D_a\der{L}{\f^\mu_a}+
D_{a b}^2 \der{L}{\f^\mu_{a b}} -\dots
\end{equation}

The map $E^{0,m}_0{\buildrel p\over\longrightarrow}E_1^{0,m}$ used in the construction
of~\eqref {varcomplex} corresponds to the projection of Lagrangians to equivalence
classes modulo $D$-coboundaries, ${\bf L}\mapsto [{\bf L}]={\bf L}+D(\O^{0,m-1})$. If
${\bf F}=F^a\,\volx_a$ is a horizontal $(m-1)$-form, then $D{\bf F}=D_aF^a\,\volx$, so in
the classical language, equivalent Lagrangians ``differ by a divergence''. On arbitrary
section $s$,
\begin{equation}\label{divergencedifference}
  {\bf L}^\prime\big\vert_s-{\bf L}\big\vert_s
  =D{\bf F}\big\vert_s=\der{}{x^a}\left(F^a\big\vert_s\right)\volx.
\end{equation}
Hence,  $\int {\bf L}^\prime$ and $\int {\bf L}$ may differ only by boundary terms, and
the variational derivative~\eqref{varderivative} is well-defined  on classes $[{\bf
L}]\in E^{0,m}_1$.

Now compare these classical considerations with the corresponding differential in the
variational complex~\eqref {varcomplex}. Consider the action of the differential
$d_1\circ p: E_0^{0.m}\longrightarrow E_1^{1.m}$, i.e., the action of the differential
$d_1$ on the equivalence class $[{\bf L}]\in E_1^{0,m}$ of a Lagrangian ${\bf
L}=L\,\volx$. According to~\eqref {vertderivaction},
\begin{equation}\label{vinvarderivone}
  d_1[{\bf L}]=[\d {\bf L}+\dots]=\left[
        \sum_{\s,\mu}\Gsm\der{L}{\fsm}\volx
                     +\ldots
                     \right]
\end{equation}
where dots denote differential forms belonging to the image of the differential
$D:\O^{1,m-1}\rightarrow \O^{1,m}$. For example,  for ${\bf F}=F^\a\,\volx\in
\O^{0,m-1}$, we have $d_1[D{\bf F}]=[\d D{\bf F}]=-[D\d {\bf F}]=0$ (compare with~\eqref
{divergencedifference}).

To find the correspondence between the image~\eqref{vinvarderivone} of  the differential
$d_1\circ p$ and the variational derivative~\eqref{varderivative}, we shall study the
content of the space $E^{1,m}_1$ in~\eqref{varcomplex}.

Consider the filtration in the Cartan ideal $C\O$ induced by the order of the multi-index
of Cartan forms. If $\O^{1,q}$ is the space of differential $(q+1)$-forms that are linear
in Cartan forms, then denote by $\O^{1,q}_{(k)}$ the subspace of $\O^{1,q}$ consisting of
forms  $\omega$ $=\sum_{|\s|\leq k} \Gsm  \osm$ where $\osm$ is a horizontal $q$-form.
One can show that if  $\omega\in \O^{1,q}_{(k)}$ for $k\geq 1$ and $D\omega=0$, then
$\omega$ is equal up to a $D$-coboundary to a form $\omega^\prime= $ $\o-D\tau\in
\O^{1,q}_{(k-1)}$. Thus every $D$-closed $\o\in \O^{1,q}_{(k)}$ is equivalent to a form
$\tilde\omega\in \O^{1,q}_{(0)}$. For example, if $\omega\in \O^{1,m}_{(k)}$,
$\omega=\Gsm \Bsm(x,[\f])\,\volx$ with $\s={a_1\dots a_k}$, then automatically $D\o=0$
and from~(\ref{horderivaction2}--\ref{vertderivaction2}) it follows that
\begin{multline*}
  \omega=\Gsm \Bsm\,\volx=\d\fsm\,\Bsm\,\volx=(\d D\f_{\s'}^{\mu})B^{a\s'}_{\mu}\,\volx_a=
  -(D\d \f_{\s'}^{\mu})B^{a\s'}_{\mu}\,\volx_a= \\
  -D\left(\d \f_{\s'}^{\mu}B^{a\s'}_{\mu}\,\volx_a\right)-
  \d \f_{\s'}^{\mu}\,D(B^{a\s'}_{\mu}\,\volx_a)=-\G_{\s'}^{\mu}(D_aB^{a\s'}_{\mu})\,\volx-
  D(\d \f_{\s'}^{\mu}B^{a\s'}_{\mu}\,\volx_a)
\end{multline*}
where $\s^\prime=a_2\dots a_k$. By iterating this  we come to the map
$\rho:\Omega^{1,m}\to \O^{1,m}_{(0)}$:
\begin{equation}\label{gotosimple}
  \r:\Gsm\,\Bsm\,\volx\mapsto \G^{\mu}(-1)^{|\s|}D_{\s}^{|\s|}B_{\mu}\,\volx,
\end{equation}
so that
\begin{equation}
  \r(\o)=\o + D(\O^{1,m-1}).
\end{equation}
\begin{prop}
The map $\rho$ defined by~\eqref{gotosimple} establishes an isomorphism between the space
$E^{1,m}_1$ and the subspace $\O^{1,m}_{(0)}\subset \O^{1,m}$.
\end{prop}
We call forms $\omega\in \O^{1,m}_{(0)}$ the \textit{canonical representatives} of
cohomological classes $[\o]\in E^{1,m}_1$. (If $\omega\in \O^{1,q}_{(0)}$, where $q<m$,
and $D\omega=0$, then one can notice directly that $\o=0$.  Thus we come to the statement
of Proposition~\ref{purealgebraic} in the case $q=1$. For $p\geq 2$ the argument is
similar).

Similar analysis of the contents of $E_1^{p,m}$ can be performed for $p>1$.

Using the homomorphism $\r$, we immediately deduce from~\eqref{vinvarderivone} that
\begin{equation}\label{nakonetsto}
  d_1[{\bf L}]=[\d {\bf
  L}]=\left[\sum_{\s,\mu}\Gsm\der{L}{\fsm}\volx\right]=[\G^\mu\fmu({\bf L})\,\volx],
\end{equation}
where $\fmu({\bf L})$ is the variational derivative~\eqref{varderivative} for the
Lagrangian ${\bf L}=L\,\volx$.

Relations between these purely algebraic considerations and the variational problem for
the functional~\eqref{actionfunctional} can be established with the help of evolutionary
vector fields considered above. For an evolutionary field
$Y=Y^\mu(x,[\f])\,\lder{}{\f^{\mu}}$ consider the  maps $I_Y:E^{p,m}_1\to E^{p-1,m}_1$
and $\lie{Y}:E^{p,m}_1\to E^{p,m}_1$,
\begin{align}
  I_Y [\o]&:=\left[\iota_{\P{Y}}\o\right],\label{secondoperation} \\
  \lie{Y}[\o]&:=\left[\lie{\P{Y}}\o\right].
\end{align}
where $\P{Y}\in\Vect(J)$ is defined by~\eqref{evolutionary2}. These maps are
well-defined, because $\iota_{\P{Y}}$ and $\lie{\P{Y}}$ commute with $D$
(see~\eqref{Dwithev} and \eqref{lieevwithD}), so for an arbitrary form $\tau\in
\O^{p,m-1}$, $\iota_{\P{Y}}D\tau=-D(\iota_{\P{Y}}\tau)$ and
$\lie{\P{Y}}D\tau=D(\lie{\P{Y}}\tau)$. From~\eqref{deltawithev} it also follows that
\begin{equation}
  [d_1,I_Y]=\lie{{Y}}
\end{equation}
on classes of forms in $E_1^{p,m}$.

In particular, let $[\o]=\left[\Gsm \Bsm\,\volx\right]$ be an arbitrary class in
$E_1^{1,m}$ and $\omega^\prime=\r(\o)=\G^\mu \widetilde B_\mu\,\volx$ be its canonical
representative ($\widetilde B_\mu=(-1)^{|\s|} D_\s^{|s|} \Bsm$). Then it follows
from~\eqref{secondoperation} that for an arbitrary evolutionary vector field $Y$
\begin{equation}\label{itisdifferential}
  I_Y[\o]=I_Y[\o']=\left[Y^{\mu}\widetilde B_\mu\,\volx\right].
\end{equation}
Consider a Lagrangian ${\bf L}=L\,\volx$. Let $Y=Y^\mu(x,[\f])\,\lder{}{\f^{\mu}}$ be an
arbitrary evolutionary  vector field. For a section $s(x)=(x,\f(x))$  of the bundle $E$,
the value of $Y$ on $s$ gives an infinitesimal variation:
\begin{equation}
  \f^{\mu}(x)\mapsto \f^{\mu}_{\e}(x)=\f^{\mu}(x)+\e Y^{\mu}(x,[\f(x)]).
\end{equation}
Then
\begin{equation}
  {\bf L}\big\vert_s\mapsto {\bf L}\big\vert_{s_{\e}}={\bf L}\big\vert_s+\e\left(\lie{\P{Y}}{\bf
  L}\right)\big\vert_s
\end{equation}
and
\begin{equation}
  S[\f]\mapsto S[\f_{\e}]=S[\f]+\e\int \left(\lie{\P{Y}}{\bf L}\right)\big\vert_s.
\end{equation}
From the results above we know that $\lie{\P{Y}}({\bf L})=\iota_{\P{Y}}\d{\bf
L}+\d\iota_{\P{Y}}{\bf L}=\iota_{\P{Y}}\d{\bf L}$. It follows
from~\eqref{itisdifferential} and \eqref{nakonetsto} that
\begin{equation}
  \left[\lie{\P{Y}}({\bf L})\right]=\left[\iota_{\P{Y}}\d{\bf L}\right]=I_Y\left[\d{\bf
  L}\right]=\left[Y^{\mu}\,\fmu({\bf L})\,\volx\right].
\end{equation}
Since the integral is constant on classes (for compactly supported forms), we conclude
that
\begin{equation}\label{genuinevariation}
  S[\f_{\e}]=S[\f]+\e\int \left(Y^{\mu}\,\fmu({\bf L})\right)\big\vert_s\,\volx.
\end{equation}
In other words, these  considerations make it possible to recover in a purely algebraic
way the Euler--Lagrange formula for the  variation of  action induced by a variation of
section.

To summarize: working with classes in $E_1^{p,m}$ is an algebraic model of integration.
Passing to other representative inside a class, e.g., getting $\G^{\mu}$'s from $\Gsm$'s
as above, corresponds to classical integration by parts. With this in mind, we rewrite
the variational complex~\eqref{varcomplex} as
\begin{equation}\label{varcomplex2}
 0\to\O^{0}_{\text{hor}}  {\buildrel D\over\longrightarrow}
   \O^{1}_{\text{hor}}  {\buildrel D\over\longrightarrow}\dots
              {\buildrel D\over\longrightarrow}
   \O^{m}_{\text{hor}}  {\buildrel \d\over\longrightarrow}
   {\textstyle\int} \O^{1,m}{\buildrel \d\over\longrightarrow}
   {\textstyle\int} \O^{2,m}{\buildrel \d\over\longrightarrow}
   {\textstyle\int} \O^{3,m}{\buildrel \d\over\longrightarrow}\dots
\end{equation}
with  a suggestive notation $\int \O^{p,m}:=E_1^{p,m}=\O^{p,m}/D(\O^{p,m-1})$ and
restoring $\d$ for $d_1$, as acting on classes. The formula~\eqref{nakonetsto} can be
rewritten then as
\begin{equation}\label{abstractvariation}
  \d[{\bf L}]=\left[\d\f^{\mu}\,\fmu({\bf L})\,\volx\right],
\end{equation}
which is identical (even to the letter $\d$) with the classical formula for the variation
of the functional $S$ (which we can identify with $[{\bf L}]$). The Cartan forms
$\G^{\mu}=\d\f^{\mu}$ play the role of ``abstract variations'', similar to differentials
$dx^a$ of the usual calculus, while evolutionary vector fields correspond to actual
variations of sections, analogous to vectors in the usual calculus. The
formula~\eqref{genuinevariation} corresponds to ``taking value'' of the
``differential''~\eqref{abstractvariation} on a ``vector'' $Y$.

\section{The complex of variational derivatives and relations of two complexes}

\def \Euler {{\rm \bf E}}

\subsection{The complex of variational derivatives}

In this section we shall consider another complex related with Euler-Lagrange equations:
the complex of Lagrangians of parametrized surfaces in a given manifold $M$. Actually, we
shall perform all considerations for an arbitrary supermanifold. The variational complex
considered above can be generalized for supercase, too (see
subsection~\ref{subsecvarsuper}).

Let $M=M^{m|n}$ be a supermanifold with coordinates $x^a$. Consider an $r|s$-dimensional
coordinate superspace  $\R{r|s}$ with standard coordinates  $t^i$. In our notation some
of coordinates are even, some odd. The parity of indices   is the parity of the
respective coordinates, $\at:=\tilde x^a$, etc. In the sequel we shall often omit the
prefix ``super'' (unless required for clarity). Consider smooth maps of a neighborhood of
zero in $\R{r|s}$ to $M$.   Consider the supermanifold of jets of such maps, of order
$k$, at point zero. (In the supercase jets are defined exactly as in the purely even
case.) Denote it by $T_{r|s}^{(k)}M$. In local coordinates the elements of
$T_{r|s}^{(k)}M$ will be $[x]=(x^a_\s)=(x^a,x^a_i,...,x^a_{i_1\dots i_k})$. There are
natural projections $T_{r|s}^{(k+1)}M\to T_{r|s}^{(k)}M$. Consider the inverse limit
$T^{\infty}_{r|s}M$ of the sequence of bundles:
\begin{equation}
  \ldots\to T_{r|s}^{(k+1)}M\to T_{r|s}^{(k)}M\ldots\rightarrow
    T_{r|s}^{(1)}M\rightarrow T_{r|s}^{(0)}M=M.
\end{equation}

\begin{de}
We call the manifold $T_{r|s}^{(k)}M$  the manifold of \textit{tangent $r|s$-elements of
order $k$} ($k=0,1,2,\ldots,\infty$).
\end{de}

Tangent elements of infinite order will be shortly called \textit{tangent elements}.
Clearly, tangent $1|0$-elements of first order are simply (even) tangent vectors:
$T_{1|0}^{(1)}M=TM$. Tangent $r|s$-elements of first order are arrays of $r$ even and $s$
odd tangent vectors. Notice that for $k>1$, the bundle $T_{r|s}^{(k)}M\to M$ is not a
vector bundle.

We define \textit{$r|s$-Lagrangians} on $M$ as  smooth functions  on the space of
$r|s$-tangent elements $T_{r|s}^\infty M$. (As before, we consider functions of finite
order.) Denote the space of all \textit{$r|s$-Lagrangians} by $\F^{r|s}=\F^{r|s}(M)$.

Consider \textit{$r|s$-paths} ($r|s$-dimensional parametrized surfaces)  $\g:U^{r|s}\to
M$ (where $U^{r|s}\subset \R{r|s}$). Every $r|s$-path $x^a=x^a(t)$ at the point $t$
defines the tangent element
\begin{equation}
  [x(t)]=\left(x^a(t),\der{x^a}{t}(t),\dder{x^a}{t^i}{t^j}(t), \ldots\right).
\end{equation}
Every Lagrangian $L\in\F^{r|s}(M)$ defines a functional on {$r|s$-paths}:
\begin{equation}\label{sotgamma}
    S[\g]=\int_{U^{r|s}} L([x(t)])\,\volt,
\end{equation}
where  $\volt$ is the standard coordinate volume form on $U^{r|s}\subset \R{r|s}$.

\begin{rem}
There are two essential differences with the setup of the previous section. First, the
space of parameters $\R{r|s}$ is endowed with fixed coordinates $t^i$. This allows to
consider the ``standard'' volume element $\volt$ and to define our Lagrangians as
functions (not as forms on jet space). Second, Lagrangians considered here are
geometrical objects on $M$ (not $\R{r|s}$ or $\R{r|s}\times M$), hence cannot depend on
$t^i$.
\end{rem}

From the variational problem for the functional~\eqref{sotgamma} one can easily obtain
the following formula for the variational derivative:
\begin{equation}\label{varderivpoputi}
  \fa(L)=\der{L}{x^a}-(-1)^{\at\itt}D_i\left(\der{L}{\xia}\right)+\ldots =
  \sum_{|\s|=0}^{\infty} (-1)^{|\s|+\at\tilde\s}D_{\s}^{|\s|}\left(\der{L}{\xia}\right).
\end{equation}
Here $D_i$ stands for the total derivative with respect to $t^i$ (covariant derivative in
the terminology of the previous sections):
\begin{equation}
  D_i=\xia\der{}{x^a}+\xija\der{}{\xja}+\ldots .
\end{equation}
Clearly, $\fa(L)=\fa(\L)$ for the Lagrangian $\L=L\,\volt$ considered as a form (see the
next subsection on the variational complex in the supercase.)

\begin{de}
The \textit{differential} of an $r|s$-Lagrangian $L$ is the  $(r+1|s)$-Lagrangian $\dd L$
defined by the formula~\cite{tv:lag}
\begin{equation}                  \label{varderivative2}
            \dd L=x^{a}_{r+1}\fa(L).
\end{equation}
\end{de}
\noindent Here $\dd L$ depends on $x^a_{i_1\dots i_k}$ where the lower indices run over
the set that includes one more even index $r+1$, corresponding to the new even variable
$t^{r+1}$. We call the operator $\dd$ the \textit{variational differential}.

\begin{prop}[see~\cite{tv:lag}]
$\dd\,^2=0$.
\end{prop}

We arrive at the \textit{complex of variational derivatives} on a (super)manifold $M$:
\begin{equation}
             \label{varderivativescomplex}
    0\to \F^{0|s}{\buildrel \dd\over\longrightarrow}\F^{1|s}
    {\buildrel \dd\over\longrightarrow}\F^{2|s}{\buildrel \dd\over\longrightarrow}\ldots
\end{equation}

In the purely even case, an example of $r$-Lagrangian of first order is provided by an
arbitrary $r$-form. It can be shown that in this case the operator $\dd$ generalizes the
usual Cartan-de~Rham exterior differential of forms. In the supercase, the correct
definition of the Cartan-de~Rham complex
(see~\cite{tv:compl},\cite{tv:git},\cite{tv:dual}) is based on this construction. We
shall return to this in subsection~\ref{subsecrelation}. One of the applications of the
complex~\eqref{varderivativescomplex} is to the study of symmetry of Lagrangians and
equations of motion~\cite{hov:double}.

\subsection{Variational complex in super context}\label{subsecvarsuper}
Much of the considerations of Sections~\ref{secprelim} and \ref{secvin} carries over to
supermanifolds without difficulties. A tricky thing is the choice of the precise class of
forms that should be used. The analog of the Cartan-de~Rham complex for supermanifolds in
its full generality is highly nontrivial (see~\cite{tv:git},\cite{tv:dual}). (This was
one of the sources of the current research, in particular, of the study of the complex of
variational derivatives. See Subsection~\ref{subsecdens}) However, for the purposes of
this paper, these complications can be circumvented, as it is shown below, if we are only
interested in forms necessary to construct the super version of the variational
complex~\eqref{varcomplex2}.

Let $E\to M$ be a fibre bundle where both $E$ and $M$ are supermanifolds. Let $\dim
M=m|n$ (dimension of the base). We continue to use the coordinates $x^a$ for $M$ and
$\f^{\mu}$ for the fibre. Some of them are even, some odd. The parity for all objects is
denoted by tilde. The jet bundles, in particular, the infinite jet bundle $J\to M$, are
defined exactly as in the purely even case. The natural coordinates for jets, $\fsm$,
have parity of the corresponding partial derivatives: $\tilde\fsm=\mut+\tilde\s$, where
$\tilde\s:=\at_1+\ldots+\at_k$ for a multi-index $\s=a_1\ldots a_k$. They enjoy the
symmetry property $\f^{\mu}_{\s ab\tau}=(-1)^{\at\bt}\f^{\mu}_{\s ba\tau}$, for each pair
of neighbor indices. As before, all functions are of finite order, i.e., depend on a
finite number of the coordinates $\fsm$. Vector fields can have infinitely many nonzero
coefficients.

In the context of this paper it is possible to simplify our task by considering a
particular class of Bernstein-Leites pseudodifferential forms instead of arbitrary super
forms. Let $\Pi TJ$ be the antitangent bundle for $J$, i.e., the tangent bundle with
reversed parity of fibres. The natural local coordinates in $\Pi TJ$ are $x^a$, $\fsm$,
$dx^a$, $d\fsm$. The differentials are treated as independent variables of parity
opposite to that of the respective coordinates. We shall consider various functions on
$\Pi TJ$ including generalized functions w.r.t. the variables $dx^a$, $d\fsm$. The space
of all such functions (without specifying a class) will be denoted here by $\O(J)$.
Warning: it is not an algebra, as not for all functions  multiplication is defined. In
particular, inside $\O(J)$ is contained the algebra of ``na\"\i ve'' differential forms,
consisting of functions polynomial in differentials. The exterior differential $d$ can be
viewed as an odd vector field on the supermanifold $\Pi TJ$, $d(x^a)=dx^a$,
$d(\fsm)=d\fsm$, $d(dx^a)=0$, $d(d\fsm)=0$.

The Cartan connection form is defined as before,
\begin{equation}\label{supercartanconnection}
  \G=\Gsm\der{}{\fsm},
\end{equation}
where
\begin{equation}\label{supercartanform}
  \Gsm=d\fsm-dx^a\fasm.
\end{equation}
It defines the decomposition of tangent and cotangent spaces and their opposites (spaces
with reversed parity) into the horizontal and vertical subspaces. In particular, for the
tangent space a basis of the horizontal subspace consists of partial covariant
derivatives
\begin{equation}\label{superpartcov}
  D_a=\der{}{x^a}+\fasm\der{}{\fsm},
\end{equation}
and a basis of the vertical subspace consists of the derivatives $\lder{}{\fsm}$. For the
anticotangent space,  a basis of the horizontal subspace consists of the differentials
$dx^a$ and a basis of the vertical subspace consists of the forms $\Gsm$. These bases of
vectors and $1$-forms are dual. It is convenient to use $dx^a$ and $\Gsm$ as fibre
coordinates in $\Pi TJ$ instead of $dx^a$ and $d\fsm$. Notice that
\begin{equation}
  d\Gsm=(-1)^{\at}dx^a\Gasm.
\end{equation}

The form-valued vector field $\G$ gives the Lie derivative $\lie{\G}$, which is a vector
field on $\Pi TJ$ generating the infinitesimal transformation
\begin{align}
  \fsm & \mapsto \fsm+\e\Gsm \\
  d\fsm & \mapsto d\fsm-\e d\Gsm=d\fsm-\e (-1)^{\at}\Gasm
\end{align}
(here $\e$ is an odd constant, $\e^2=0$). It follows that $\Gsm\mapsto \Gsm$, i.e.,
\begin{equation}
  \lie{\G}\Gsm=0,
\end{equation}
which is equivalent to the vanishing of the Nijenhuis bracket
\begin{equation}
  [\G,\G]=0
\end{equation}
(flatness of the connection $\G$). Hence we can decompose $d=D+\d$:
if $\o=\o(x,[\f],dx,\G)$ is a function on $\Pi TJ$ written in the coordinates $x^a$,
$\fsm$, $dx^a$, $\Gsm$, then
\begin{equation}\label{decompofd}
  d\o=\underbrace{dx^a\left(\der{\o}{x^a}+\fasm\der{\o}{\fsm}+(-1)^{\at}\Gasm\der{\o}{\Gsm}\right)}_{\text{horizontal}}
+\underbrace{\Gsm\der{\o}{\fsm}}_{\text{vertical}},
\end{equation}
or
\begin{align}
  D & = dx^a\left(\der{}{x^a}+\fasm\der{}{\fsm}+(-1)^{\at}\Gasm\der{}{\Gsm}\right), \label{superD}\\
  \d &=\lie{\G} =\Gsm\der{}{\fsm}. \label{superdelta}
\end{align}
The differentials $D$ and $\d$ commute: $D\d+\d D=0$, and $D^2=\d^2=0$. In particular,
\begin{equation}
  \Gsm=\d\fsm.
\end{equation}

One can introduce a bigrading into the algebra of na{\"\i}ve differential forms exactly
as in the purely even case considered in Section~\ref{secprelim}:
\begin{equation}\label{superbigrading1}
  \O^{p,q}(J):=\Bigl\{\o\in\O(J)\,\bigl\vert\, \#\Gsm(\o)=p \text{  \,and\,  } \#dx^a(\o)=q
  \bigr.\Bigr\},
\end{equation}
where $\o$ is assumed to be polynomial in $dx^a$ and $\Gsm$, and $\#dx^a$, $\#\Gsm$ stand
for the degree in the corresponding variables. It is possible to carry on with the
filtration and spectral sequence as above, but the trouble is that in the super case
these considerations give nothing for the variational problem: since na{\"\i}ve
differential forms cannot be integrated over supermanifolds
(see~\cite{tv:git},\cite{tv:dual} for discussion), Lagrangians are not contained in the
complex $\O^{**}(J)$.

The correct class of forms on $J$, adequate for our purposes, can be introduced as
follows. (We are going to introduce a new class of forms particularly designed for the
needs of variational complex. They will be some hybrids of integral forms and na\"{\i}ve
differential forms, see remark below.)

Recall that for arbitrary ``super'' variables $z^a$ the delta-function $\d(z)$
corresponds to the distribution $f\mapsto f(0)$ w.r.t. the Berezin integration in
variables $z^a$. In particular, for even variables it is the usual delta-function, and
for odd variables $\x^i$
\begin{equation}\label{deltaodd}
  \d(\x^1,\ldots,\x^n)=\x^n\x^{n-1}\ldots\x^1
\end{equation}
(maximal product; notice the inverse order). Delta-function $\d(z)$ satisfies the
property
\begin{equation}\label{deltazamena}
  \d(z)=\d(z')\cdot\Bigl(\Ber\der{z}{z'}\Bigr)^{\!-1}\sign\det\Bigl(\der{z}{z'}\Bigr)_{\!00},
\end{equation}
for a non-singular change of variables, where $\Ber$ is the Berezinian and the indices
$00$ refer to the even-even block of the Jacobi matrix. Consider the delta-function
$\d(dx)$ of the variables $dx^a$ (there should be no confusion with the vertical
differential $\d$!). This is well-defined, because the differentials $dx^a$ transform
through themselves under changes of variables $x^a$, $\fsm$. Moreover, from the
property~\eqref{deltazamena} we obtain the following transformation law:
\begin{multline}\label{deltadxzamena}
  \d(dx)=\d(dx')\cdot\Bigl(\Ber\Bigl(\der{x}{x'}\Bigr)^{\!\!\Pi}
  \Bigr)^{\!-1}\sign\det\left(\der{x}{x'}\right)_{\!\!11}=\\
  \d(dx')\cdot\Ber\der{x}{x'}\cdot\sign\det\left(\der{x}{x'}\right)_{\!\!11},
\end{multline}
where the superscript $\Pi$ denotes reversion of parity of rows and columns of a matrix
(notice that $\Ber g^{\Pi}=(\Ber g)^{-1}$). It follows that up to a sign factor, $\d(dx)$
transforms exactly as the Berezin volume element $\volx$ and thus can be identified with
such. More precisely, we consider (generalized) functions on $\Pi TJ$ that are supported
at the closed submanifold $\Pi VJ\hookrightarrow\Pi TJ$ ($VJ$ is the vertical subbundle),
which is  locally specified by the equations $dx^a=0$. Let them take values in the local
system $\sign\det (TM)_1$. This eliminates the above sign factor. Now, every such
function is a linear combination of the delta-function $\d(dx)$ and its derivatives. We
introduce the following notation:
\begin{align}
  \volx & :=\d(dx) \label{volx}\\
  \volx_{a_1\dots a_k} & :=\der{}{dx^{a_1}}\ldots \der{}{dx^{a_k}}\,\d(dx).\label{volxa}
\end{align}
After introducing the said local coefficients, the formula~\eqref{volx} is a genuine
identification. (In particular, for a purely even manifold $M^m$, we have
$\volx=\d(dx)=dx^mdx^{m-1}\dots dx^1=\pm dx^1\dots dx^m$; notice slight difference of
sign from our notation in Section~\ref{secvin}.) The forms that we need have the
appearance:
\begin{equation}\label{nashiformy}
  \o=\frac{1}{k!}\,\o^{a_1\ldots a_k}\volx_{a_k\ldots a_1}
\end{equation}
($k!$ and the order of indices are chosen for convenience). The coefficients
$\o^{a_1\ldots a_k}$ depend on $x^a$, $\fsm$, $\Gsm$. From~\eqref{volx} and \eqref{volxa}
follow the rules of operations with the symbols $\volx$ and $\volx_{a_1\dots a_k}$:
\begin{align}
  dx^a\,\volx & =0 \label{rules1} \\
  \der{}{dx^a}\,\volx_{a_1\dots a_k} & =\volx_{aa_1\dots a_k} \label{rules2}\\
  dx^a\,\volx_{a_1\dots a_k}&=\sum_{i=1}^{k} (-1)^{\at+(\at+1)(\a_1+\ldots+\at_{i-1}+i-1)}
  \d_{a_i}^a\,\volx_{a_1\dots \hat{a_i}\ldots a_k}.\label{rules3}
\end{align}
where hat means that the index is dropped. In particular, for arbitrary $a$,
\begin{equation}\label{dxavolxa}
  \underbrace{dx^a\volx_a}_{\text{no summation}}=(-1)^{\at}\volx.
\end{equation}
The space of forms~\eqref{nashiformy} is not an algebra, but as follows
from~\eqref{rules1}--\eqref{rules3}, it is a module over the ``na\"\i ve'' algebra of
polynomial forms, of ``Fock type'', where multiplication by $dx^a$ acts as annihilation
operators for a ``vacuum vector'' $\volx$.

Consider forms~\eqref{nashiformy} with coefficients polynomial in $\Gsm$. Define
\begin{equation}\label{superbigrading2}
  \S^{p,q}(J):=\Bigl\{\o=\frac{1}{k!}\,\o^{a_1\ldots a_k}\volx_{a_k\ldots a_1}
  \,\bigl\vert\, \#\Gsm(\o)=p \text{  \,and\,  } k=m-q
  \bigr.\Bigr\}.
\end{equation}
(We have $\S^{p,q}(J)\subset\O(J)\otimes{\cal E}$, where ${\cal E}$ is the orienting
sheaf for $(TM)_1.$) Here $m$ is the even dimension of the base: $M=M^{m|n}$. Operations
$D$ and $\d$ act on $\S^{**}$. For $\o$ as in~\eqref{nashiformy}, one can find
from~\eqref{superD}--\eqref{superdelta} and \eqref{rules1}--\eqref{rules3}  that $D$ acts
as a ``divergence'':
\begin{equation}
  D\o  =-\frac{1}{(k-1)!}(-1)^{(\ot+m)(\at+1)}D_a\o^{ab_1\dots b_{k-1}}\,\volx_{b_{k-1}\dots
  b_1},
\end{equation}
and $\d$ acts componentwise:
\begin{equation}
  \d\o=\frac{1}{k!}\,\Gsm\der{\o^{a_1\ldots a_k}}{\fsm}\,\volx_{a_k\ldots a_1}.
\end{equation}
Notice that $D:\S^{p,q}\to\S^{p,q+1}$ and $\d:\S^{p,q}\to\S^{p+1,q}$.  The bicomplex
$\S^{**}$ is bounded at the bottom and at the right: $\S^{p,q}=0$ for $p<0$ or $q>m$. In
particular, there is a complex of horizontal forms
\begin{equation}\label{superhorizontalcomplex}
 \dots\to \S_{\rm hor}^0{\buildrel D\over\longrightarrow}\,
     \S_{\rm hor}^1{\buildrel D\over\longrightarrow}\,
                   \dots
        {\buildrel D\over\longrightarrow}\,
        \S_{\rm hor}^{m-1}{\buildrel D\over\longrightarrow}\,
                  \S_{\rm hor}^m\to 0,
\end{equation}
where $\S^{q}_{\text{hor}}=\S^{0,q}$. (It is not bounded at the left.)

Applying the machinery of homological algebra to the bicomplex $\S^{**}(J)$, we arrive at
a spectral sequence analogous to~\eqref{specsequence}. In particular, we obtain the
\textit{variational complex} for the super case:
\begin{equation}\label{supervarcomplex2}
 \dots\to\S^{0}_{\text{hor}}  {\buildrel D\over\longrightarrow}
   \S^{1}_{\text{hor}}  {\buildrel D\over\longrightarrow}\dots
              {\buildrel D\over\longrightarrow}
   \S^{m}_{\text{hor}}  {\buildrel \d\over\longrightarrow}
   {\textstyle\int} \S^{1,m}{\buildrel \d\over\longrightarrow}
   {\textstyle\int} \S^{2,m}{\buildrel \d\over\longrightarrow}
   {\textstyle\int} \S^{3,m}{\buildrel \d\over\longrightarrow}\dots .
\end{equation}
Here the elements of the space ${\textstyle\int} \S^{p,m}$ are formal integrals $\int
\o$, where $\o\in\S^{p,m}$.  These formal integrals are in $1-1$-correspondence with
classes $\o\mod D(\S^{p,m-1})$. One subtlety is that for consistence with the parity of
the genuine Berezin integral we define the isomorphism ${\textstyle\int} \S^{p,m}\cong
E^{p,m}_1=\S^{p,m}/D(\S^{p,m-1})$ as having parity $n$ (where $n$ is the odd dimension of
the base $M=M^{m|n}$).  Notice the following properties of this complex: it is not
bounded from the left (differently from~\eqref{varcomplex2}) and it is not bounded from
the right (similar to~\eqref{varcomplex2}). Still, its cohomology is the ordinary
cohomology of the underlying topological space of the bundle $E$.

The relation with the classical calculus of variations is exactly as in the purely even
case. Suppose we have a Lagrangian $\L=L(x,[\f])\,\volx$. It is an element of $\S^{0,m}$.
The action (treated formally) is the class $\int \L\in \int\S^{0,m}\cong
\S^{0,m}/D(\S^{0,m-1})$. Consider $\d\int \L$ (since $\d$ and $D$ commute, it makes sense
to apply $\d$  modulo $D$-coboundaries). We have:
\begin{equation}\label{supervar1}
  \d\int\L=(-1)^{n}\int \d\L=(-1)^{n}\int \Gsm\der{L}{\fsm}\,\volx.
\end{equation}
To transform this, consider first an arbitrary form $\o\in\S^{1,m}$. Let
\begin{equation}
  \o=\G^{\mu}_{a_1\dots a_k} B^{a_1\dots a_k}_{\mu}\,\volx=
  \d\f^{\mu}_{a_1\dots a_k} B^{a_1\dots a_k}_{\mu}\,\volx.
\end{equation}
Exactly as we did to obtain~\eqref{gotosimple}, we can show, using~\eqref{dxavolxa}, that
there is an equality modulo $D(\S^{1,m-1})$:
\begin{equation*}
  \omega=\d\f^{\mu}_{a_1\dots a_k} B^{a_1\dots a_k}_{\mu}\,\volx
  =-(-1)^{\at_k\mut}\d\f^{\mu}_{a_1\dots a_{k-1}} D_{a_k}B^{a_1\dots
  a_k}_{\mu}\,\volx.
\end{equation*}
From here we deduce the map $\r$ giving canonical representatives: $\o\equiv\r(\o)\mod
D(\S^{1,m-1})$, where
\begin{equation}\label{supergotosimple}
  \r(\o)=(-1)^{|\s|+\tilde\s\mut}\d\f^{\mu}\,D_{\s}^{|\s|}B^{\s}_{\mu}\,\volx,
\end{equation}
where $\s=a_1\dots a_k$. Applying this to~\eqref{supervar1}, we obtain
\begin{equation}\label{supervar2}
  \d\int\L=(-1)^{n}\int
  \d\f^{\mu}(-1)^{|\s|+\tilde\s\mut}D_{\s}^{|\s|}\der{L}{\fsm}\,\volx,
\end{equation}
which is exactly
\begin{equation}
  \d\int \L=(-1)^n\int \d\f^{\mu} \fmu (\L)\,\volx,
\end{equation}
with
\begin{equation}
  \fmu (\L)=(-1)^{|\s|+\tilde\s\mut}D_{\s}^{|\s|}\der{L}{\fsm}
\end{equation}
being the variational derivative in the super case.

\begin{rem}
Forms that we introduced here are  hybrids of Bernstein--Leites integral forms
(multivector densities of weight $1$) on the base $M$ with  differential forms on fibres.
This can be seen directly if one applies the Fourier-Hodge transform to our
formulas~\eqref{volx}--\eqref{nashiformy}. We preferred not to work with a ``hybrid''
definition explicitly in order to avoid a ``nonlocal'' transformation law under
coordinate changes (though it would not occur if one sticks to the non-holonomic frame
$\Gsm$), choosing instead the language of generalized functions. To put this into a
proper framework, one should notice that, in general, the space of forms on a
supermanifold is $\Om{r}{s}$ (see~\cite{tv:dual}), where $r|s$ is a super dimension. For
the jet space $J$ the most general bicomplex has the bigrading like
$\boldsymbol{\O}^{r|s,p|q}$. The na\"{\i}ve spaces of forms $\O^{p,q}$ that we initially
introduced are exactly $\boldsymbol{\O}^{p|0,q|0}$; the ``correct'' space $\S^{p,q}$ that
we used for the variational complex coincides with $\boldsymbol{\O}^{p|0,q|n}$ (here $n$
is the odd dimension of the base $M^{m|n}$).
\end{rem}

\subsection{Relation of complexes}
\label{subsecrelation}

Let us compare the variational complex~\eqref{varcomplex}, \eqref{varcomplex2} with the
complex of variational derivatives~\eqref{varderivativescomplex}. For simplicity consider
a purely even case. (The general super case is similar.)

Consider for every $r$ the trivial fiber bundle $\pi_{(r)}=\mathbb R^r\times M$ with base
$\mathbb R^r$ and fibre  $M$ and the corresponding space of jets $J^\infty(\pi_{(r)})$.
In the sequel we shall shortly denote $J_{(r)}:=J^\infty(\pi_{(r)})$. There is a natural
projection $J_{(r)}\to T^{(\infty)}_r$ onto the bundle of tangent elements of $M$.

Consider the embedding $\O^*_{(r)}\hookrightarrow$ ${\O_{(r+1)}^*}$ induced by the
natural projection $p_r: \pi_{(r+1)}\to\pi_{(r)}$, where ${\O_{(r)}^*}$ is the space of
differential forms on the space   $J_{(r)}$. Consider the filtrations
\eqref{vinfiltration} generated by the Cartan ideals in the spaces  $\O^*_{(r)}$ and
$\O^*_{(r+1)}$ and the zero terms $E^{**}_{(r)0}$, $E^{**}_{(r+1)0}$ of the corresponding
spectral sequences (see~\eqref{vinzeroterm}).

\begin{prop}
1. Under the projection $p_r$, the $k$-th power of the Cartan ideal in $\O^*_{(r)}$ maps
to the $(k-1)$-th power of the Cartan ideal in the space $\O^*_{(r+1)}$:
\begin{equation}
      \label{mapofcartanidelas}
      p_r^*: C^k\O^*_{(r)}\longrightarrow C^{k-1}\O^*_{(r+1)}.
\end{equation}

2. The map~\eqref{mapofcartanidelas} induces a map of the zeroth terms of the
corresponding spectral sequences
                $E^{p,q}_{(r)0}
         \rightarrow E^{p-1,q+1}_{(r+1)0}$.
Thus it defines a map of bicomplexes $\O^{**}_{(r)}\to\O^{**}_{(r+1)}$ of bidegree
$(-1,1)$:
\begin{equation}
             \label{mapofbicomplexes}
  \kappa_r:\O^{p,q}_{(r)}\rightarrow
                \O^{p-1,q+1}_{(r+1)}.
\end{equation}
\end{prop}
\begin{proof}
Let $\o=\o(x,dt,\G)$ be a form in $\O^{p,q}(J_{(r)})$. Here $\G=\G_{(r)}$ stands for
Cartan forms on $J_{(r)}$. Notice that
$\G_{(r)\s}^a=\G_{(r+1)\s}^a+dt^{r+1}x^a_{r+1,\s}$. Hence
\begin{equation}
  p_r^*\o=\o(x,dt,\G_{(r)})=
  \o(x,dt,\G_{(r+1)})+dt^{r+1}x^a_{r+1,\s}{\p{\o}\over\p{\G^a_{(r)\s}}}.
\end{equation}
The first term in the r.h.s. belongs to $\O^{p,q}(J_{(r+1)})$ and the second term in the
r.h.s. belongs to $\O^{p-1,q+1}(J_{(r+1)})$. It follows that the induced map
$\kappa_r:\O^{p,q}(J_{(r)})\to\O^{p-1,q+1}(J_{(r+1)})$ is given by the formula
\begin{equation}
  \kappa_r\o=dt^{r+1}x^a_{r+1,\s}{\p{\o}\over\p{\G^a_{(r)\s}}}.
\end{equation}
\end{proof}

For example, let $\omega=\G_{(r)}^a\G^b_{(r)i}\in\O^{2,0}_{(r)}$ be a differential form
on $J_{(r)}$, where
\begin{equation}
        \label{cartanformsindifferentspaces}
             \G^a_{(r)}=
                dx^a-
          \sum_{i\leq r} dt^i x^a_{i}\,,\quad
     \G^b_{(r)i}=dx^b_i-\sum_{j\leq r} dt^jx^b_{ij}.
\end{equation}
Under the embedding \eqref{mapofcartanidelas} the form $\omega=\G^a_{(r)}\G^b_{(r)i}$
gives the form
\begin{equation}
\label{cartanformsindifferentspacestwo}
                   p^*_r\o= \left(
        \G^a_{(r+1)}+dt^{r+1} x^a_{r+1}
                    \right)
                     \left(
         \G^b_{(r+1)i}+dt^{r+1} x^b_{r+1,i}
                    \right)\in C^1\O^2_{(r+1)}
\end{equation}
on the space $J_{(r+1)}$.  Hence, $\kappa_r\o=x^a_{r+1}dt^{r+1}\G^b_{(r+1)i}+
\G^a_{(r+1)}x^b_{i.r+1}dt^{r+1}$, which belongs to $\O^{1,1}_{(r)}$.

Consider the induced action of $\kappa_r$ on variational complexes~\eqref{varcomplex2}.
Clearly, $\kappa_r$ vanishes on $\O_{(r)}^{0,q}$ and
$\kappa_r\bigl(\int\O_{(r)}^{1,r}\bigr)\subset \int \O_{(r+1)}^{0,r+1}$. Recall the map
$\r$ that takes classes in $\int\O_{(r)}^{1,r}$ to their canonical representatives,
see~\eqref{gotosimple}. Define the map $\chi_r:
\int\O_{(r+1)}^{1,r}\to\O_{(r+1)}^{0,r+1}$ as the composition $\kappa_r\circ\r$. We shall
use the map $\chi_r$ to analyze the relation between the variational complex and the
complex of variational derivatives~\eqref{varderivativescomplex}.

Assign to every $r$-Lagrangian $L$ on $M^m$ the horizontal $r$-form ${\bf
L}=L\,\volt_{(r)}$ on $J_{(r)}$, where $\volt_{(r)}$ is the coordinate volume form on the
space of parameters $\mathbb R^r$. Consider the $(r+1)$-Lagrangian $\dd L=x^a_{r+1}{\cal
F}_a(L)$, where ${\cal F}_a({L})=\fa(\L)$ is the variational
derivative~\eqref{varderivpoputi}. To it corresponds the form $\dd
L\,\volt_{(r+1)}=x^a_{r+1}{\cal F}_a(L)\,\volt_{(r+1)}$ belonging to
$\O_{(r+1)}^{0,r+1}$. On the other hand, consider the canonical representative of the
element $\d[{\bf L}]\in E^{1,r}_{(r)1}$, which is the form $\rho(\d[{\bf L}])$ $={\G^a
\cal F}_a \volt_{(r)} =(dx^a-\sum_{1\leq i\leq r}dt^i x^a_i) {\cal F}_a \,\volt_{(r)}$
belonging to $\O_{(r)}^{1,r}$. Applying $\kappa_r$ we obtain the action of $\chi_r$ on
$\d[L\,\volt_{(r)}]$. Thus we have the following proposition:

\begin{prop} For every $r$-Lagrangian $L$
\begin{equation}\label{samajavazhnajaformula}
  \dd L\,\volt_{(r+1)}=\chi_r\left(\d[L\,\volt_{(r)}]\right).
\end{equation}
\end{prop}

\begin{rem}
Vanishing of the induced map $\kappa_r:\int\O_{(r)}^{1,r}\to\int\O_{(r+1)}^{0,r}$ on the
image of the map $\d: \O_{(r)}^{0,r}\to\int\O_{(r)}^{1,r}$ together with the
formula~\eqref{samajavazhnajaformula} implies that $\dd L\,\volt_{(r+1)}$ is a
$D$-coboundary, i.e., $\dd L$ is a total divergence.
\end{rem}

\subsection{The complex of variational derivatives and covariant Lagrangians}

\label{subsecdens}

Consider $r|s$-Lagrangians $L$  on a supermanifold  $M$ such that for every $r|s$-path
$\g:x^a=x^a(t)$ the corresponding functional $S[\g]$ depends only on the image of the
path $x(t)$.  That means that if $\g_1$ and $\g_2$ are two arbitrary paths such that
$x_2(t)=x_1(f(t))$, where $f$ is an orientation preserving diffeomorphism of $\mathbb
R^{r|s}$, then
\begin{equation}
         \label{conditionondensitiesintegral}
         \int L\big\vert_{\g_1}\volt=
         \int L\big\vert_{\g_2}\volt.
\end{equation}
\begin{de} \label{0}
We call  $L$ a \textit{covariant $r|s$-Lagrangian of weight $\rho$} if for arbitrary
$r|s$-path $\g$ and an arbitrary orientation preserving diffeomorphism $f$ of $\mathbb
R^{r|s}$ the following condition holds:
\begin{equation}
       \label{conditionondensities}
              L\big\vert_{f^*\g}
                    (\volt)^\rho=
                       f^*
                    \left(
              L\big\vert_{\g}
                      (\volt)^\rho
                    \right)\,.
                \end{equation}
\end{de}

If a covariant $r|s$-Lagrangian has weight $\rho=1$ then this condition  is equivalent to
condition \eqref{conditionondensitiesintegral}. A covariant $r|s$-Lagrangian $L$ of
weight $\rho$ on $M$ defines the density $L\vert_{x(t)}(\volt)^\rho$ of weight $\r$  on
an arbitrary $r|s$-dimensional surface embedded in $M$. In particular, a covariant
$r|s$-Lagrangian $L$ of weight $\rho=0$ can be considered as a function on jets of
$r|s$-surfaces in $M$.

\begin{rem}
Covariant Lagrangians play important role in quantum field theory.  The application of
covariant Lagrangians of embedded surfaces   has been pioneered in the works of
A.S.~Schwarz. See, for example,~\cite{ass:equal}.
\end{rem}

In the language of tangent elements the covariance condition~\eqref{conditionondensities}
can be stated as follows.  In the space of tangent $r|s$-elements at every point of $M$
acts the (super)group of jets of diffeomorphisms of the neighborhood of the origin in
$\R{r|s}$ that fix the origin. Denote it $G(r|s)$. The usual general linear supergroup
$GL(r|s)$ is the factor group of $G(r|s)$ w.r.t. the normal subgroup of jets of
diffeomorphisms preserving the first infinitesimal neighborhood of the origin. It makes
sense to speak about the Berezinian for such jets of diffeomorphisms, which is simply the
Berezinian of the corresponding Jacobi matrix (i.e, the pull-back of the Berezinian on
$GL(r|s)$). The covariance condition for an $r|s$-Lagrangian $L$ means that
\begin{equation}\label{covariance}
  L(g[x])=(\Ber g)^{\r} \cdot L([x]),
\end{equation}
for all $g\in G(r|s)$. This is a nonlinear analog of the covariance condition for  first
order Lagrangians~\cite{tv:git}.

The infinitesimal version of the condition~\eqref{covariance} can be written as
\begin{equation}
    \label{infinitesimalconditionondensities}
                 (-1)^{\tilde\s\Kt}D_\s^{|\s|}
                     \left(
                 K^i(t)x^a_i
                     \right)
             {\p L\over\p x^a_\s}
                       =
                    \rho\,
             (-1)^{\itt(\Kt+1)}{\p K^i(t)\over\p t^i}\,L
\end{equation}
at $t=0$, where $K=K^i(t)\,\lder{}{t^i}$ is an arbitrary  vector field vanishing at the
origin. (Clearly, both sides depend only on its jet at the point $0$.) Explicitly, the
property~\eqref{infinitesimalconditionondensities} is equivalent to the following
sequence of identities:
\begin{align}
                    x^a_i{\p L\over\p x^a_j}+
                   2x^a_{ik}{\p L\over\p x^a_{jk}}+
                   \dots&=\rho\,(-1)^{\itt}\delta^j_i L,
                   \label{densityconditions1}
                   \\
          x^a_i{\p L\over\p x^a_{jk}}+3x^a_{il}{\p L\over\p x^a_{jkl}}+
                   \dots&=0,
                   \intertext{\hfill\dots\dots\dots\dots\qquad\qquad\qquad\qquad\qquad\qquad\qquad}
  x^a_i{\p L\over\p x^a_{j_1\dots j_N}}&=0,
\end{align}
where $L$ is a covariant $r|s$-Lagrangian of weight $\rho$ and order $N$.
\smallskip

If $K=K^i(t)\,\lder{}{t^i}$ is a genuine vector field on the space of parameters $\mathbb
R^{r|s}$, not necessarily vanishing at $0$, then it follows
from~\eqref{infinitesimalconditionondensities} that for an arbitrary covariant
$r|s$-Lagrangian $L$ of weight $\rho$
\begin{equation}
    \label{infinitesimalconditionondensitiestwo}
                 K^{(\infty)}L
                    +\rho\,
             (-1)^{\itt(\Kt+1)}{\p K^i\over\p t^i}L=0\,,
\end{equation}
where $K^{(\infty)}=K^i(t)D_i-{\cal P}_{K^i x_i^a\lder{}{x^a}}$ is the prolongation  of
the vector field $K$ (see subsection~\ref{subsecevol}).

Let us consider relations between covariant Lagrangians and differential forms.

In the case of ordinary purely even manifold  $M^m$  (odd dimension is equal to zero)
one can assign to every differential $r$-form
\begin{equation}
  \omega=\frac{1}{r!}\,\omega_{a_1\dots a_r}dx^{a_1}\dots dx^{a_r}
\end{equation}
the following covariant $r$-Lagrangian (i.e., $r|0$-Lagrangian) of first order and of
weight $\rho=1$:
\begin{equation}
              \label{diffformdensitycorrespondence}
     L_\omega=\omega_{a_1\dots a_r}
               x^{a_1}_{1}\dots x^{a_r}_{r}.
\end{equation}
For every $r$-path $\g$,  $\int_{\g}L_\omega=\int_{\g}\omega$
       and
             \begin{equation}
              \label{diffformdensitycorrespondence2}
                 L_{d\omega}=\dd L_{\o}\,,
                   \end{equation}
where $d\omega$ is the usual exterior differential of the form $\omega$. Notice that
under the identification of tangent elements with jets in $J^\infty(\pi_{(r)})$, the
projection of the pull-back of an $r$-form $\omega$ to on $J^\infty(\pi_{(r)})$ onto the
subspace $\O^{0,r}$ of horizontal $r$-forms is equal to the form $L_\omega (x)dt^1\dots
dt^r$.

In general, the variational  differential $\dd$ does not take covariant Lagrangians to
covariant Lagrangians and it increases the order: if $L$ has order $k$, then, generally,
the order of $\dd L$ is equal to $2k$.

\begin{prop}\label{beginningofvoronov}
If a Lagrangian $L$ on a purely even manifold $M$ is of first order and the Lagrangian
$\dd L$ is of first order too, then the Lagrangian $L$ is up to a constant a covariant
Lagrangian corresponding to a differential form: $L=L_{\omega}+c$. The Lagrangian $\dd L$
corresponds to the differential form $d\omega$.
\end{prop}

This proposition can be easily checked by straightforward calculations. It has the
following geometrical meaning. In the space of paths on the manifold $M$ consider a
collection of topologies ${\cal T}_k$  ($k=0,1,\dots$) such that sequence of paths $\g_n$
tends to $\g$ in the topology ${\cal T}_k$ if $|x_n(t)-x(t)|\to 0$ and a
\begin{equation}
  \left|{\p^kx^a_n(t)\over \p t^{i_1}\dots\p t^{i_l}}-
  {\p^kx^a(t)\over \p t^{i_1}\dots\p t^{i_l}}\right|\to 0
\end{equation}
for all  derivatives of order $l\leq k$. For a given Lagrangian $L$ of order $k$ the
corresponding functional on paths is continuous in topology ${\cal T}_k$. On the other
hand, if a first order covariant Lagrangian $L$ corresponds to a differential form:
$L=L_\omega$, then due to the Stokes theorem the functional $S[\g]$ corresponding to this
Lagrangian is continuous not only in the topology ${\cal T}_1$ but in the weaker topology
${\cal T}_0$. The converse implication leads to Proposition~\ref{beginningofvoronov} (see
\cite{ass:equal},\cite{hov:mult},\cite{hov:gayduk}).  In the supercase these
considerations yield the definition of forms~\cite{tv:compl}:
\begin{de}\label{forms}
A covariant $r|s$-Lagrangian $L$ of  first order on a supermanifold $M$ is called an
\textit{$r|s$-form} if the Lagrangian $\dd\,L$ is of first order too.
\end{de}

The Lagrangian $\dd L$ is $(r+1,s)$-form if $L$ is an $r|s$-form and the Stokes theorem
is valid. As follows from Proposition~\ref{beginningofvoronov}, on ordinary manifolds
$r|0$-forms are in $1-1$-correspondence with the usual differential $r$-forms and $\dd$
corresponds to the usual $d$. See for details~\cite{tv:git}.

\begin{rem}
The complete theory of forms in super case includes also objects defined similarly to
Definition~\ref{forms}, but in the dual setting, using so-called
copaths~\cite{tv:dual},\cite{tv:cartan1},\cite{tv:cartan2}.
\end{rem}

Consider now the action of the variational differential $\dd$ for covariant Lagrangians
of higher order.

\begin{prop} \label{trivial}
Covariant Lagrangians that are coboundaries in the complex  of variational derivatives
must be of first order.
\end{prop}
\begin{proof}
Let a coboundary $\dd\,L= x^a_{r}{\cal F}_a\in\F^{r|s}$ be a covariant $r|s$-Lagrangian.
Then it must have weight $\rho=1$, because it is linear in variables $x^a_{r}$. Let us
show  that $\dd\,L$ is of first order. Consider the identity~\eqref{densityconditions1}
for the case when index $i$ is equal to $r$. The variational derivative ${\cal F}_a$ does
not depend on variables $x^a_{i_1\dots i_k}$ if at least one of indices ${i_1,\dots,
i_k}$ is equal to $r$ and $k\geq 2$. Hence it follows from condition
\eqref{densityconditions1} that $\dd\,L$ does not depend on variables $x^a_{i_1\dots
i_k}$ if $k\geq 2$, thus, the Lagrangian $\dd\,L$ is of first order.
\end{proof}

It follows from this Proposition that a  covariant closed $r|s$-Lagrangian  $L$ of order
higher than one is a non-trivial cocycle of the complex of variational derivatives
$\F^{*|s}$.

\smallskip
Finally, we consider some constructions for higher order covariant Lagrangians.

\smallskip
\noindent{1. \textit{Lie derivative and the variational derivative of covariant
Lagrangians.}

\smallskip
\noindent Let $X=X^a(x){\lder{}{x^a}}$ be an arbitrary vector field on $M$ and $L$ be a
covariant $r|s$-Lagrangian of arbitrary weight $\rho$. Then the Lie derivative of this
Lagrangian along the vector field $X$, $\lie{X} L={\cal P}_X L=(-1)^{\tilde\s\at}
D_\s^{|\s|}X^a\lder{L}{x^a_{\s}}$, obviously,  is also a covariant $r|s$-Lagrangian of
the same weight and of the same order.

The Lagrangian $L$ and the vector field $X=X^a(x){\lder{}{x^a}}$ yield also another
$r|s$-Lagrangian $X^a(x){\cal F}_a(L)$, where  ${\cal F}(L)$ is the variational
derivative of the Lagrangian $L$ (see \eqref{varderivpoputi}). The Lagrangians ${\cal
L}_X L$ and $X^a(x){\cal F}_a(L)$ differ by a divergence: ${\cal L}_X L-X^a(x){\cal
F}_a(L)=D_iB^i$ (compare the end of subsection~\ref{subsecvinclassicrelation}). If the
covariant Lagrangian $L$ has weight $\rho=1$, then one can prove that $X^a(x){\cal
F}_a(L)$ is also covariant of weight $\rho=1$. In general, it has order $2k$, if the
Lagrangian $L$ has order $k$.

\smallskip
\noindent {2. \textit{Composition of Lagrangians.}

\smallskip
\noindent Let $L$ be an $r|s$-Lagrangian on a manifold $N$ and let $F$ be an
$r|s$-Lagrangian on a manifold $M$ taking values in the manifold $N$. Then one can
consider the $r|s$-Lagrangian $L\circ F$ called the \textit{composition} of these
Lagrangians. If $x^a$ and $y^\mu$ are local coordinates on $M$ and $N$ respectively,
$L=L([y])$  and $F:y^{\mu}=F^\mu([x])$, then consider the formal substitution
\begin{align}\label{hovikcompos}
  y^\mu & =F^\mu([x]) \\
  y^\mu_i & =D_i F^\mu([x]) \\
  y^\mu_{ij} & =D^2_{ij} F^\mu([x])
\end{align}
{\hfil\dots\dots\dots\qquad}

\noindent The Lagrangian $L\circ F$ is obtained from $L$ by this substitution.

One can see that if the Lagrangian $L$ is closed, then the Lagrangian $L\circ F$ is also
closed:
\begin{equation}\label{closednessofcomposition}
  \dd L=0\Rightarrow \dd (L\circ F)=0.
\end{equation}
Notice the special case when $F$ is a covariant $r|s$-Lagrangian of the weight $\rho$=0.
($F$ can be viewed as an $N$-valued function on jets of $r|s$-surfaces in $M$). In this
case, if $L$ is a covariant (closed) $r|s$-Lagrangian of weight $\rho$, then $L\circ F$
is also covariant (closed) Lagrangian of the same weight $\rho$. The order of Lagrangian
$L\circ F$ is equal, in general, to the sum of orders of the Lagrangians $L$ and $F$.

Consider the following toy examples of this construction.

Let $M$ be the Euclidean space $\mathbb R^m$. For every $r$, $1\leq r\leq m$ consider as
the target space $N$ the manifold of oriented $r$-dimensional linear subspaces of
$\mathbb R^m$ (the oriented Grassmannian $G^+_r=G^+_r(\mathbb R^m)$). Consider a function
$F_r$ with values in $G^+_r$ such that $F_r$ assigns to every point $x\in \mathbb R^m$
and an arbitrary oriented $r$-dimensional plane $\Pi^r$ through this point the oriented
linear subspace parallel to  $\Pi^r$. The function $F_r$ defines on $\mathbb R^m$ a
covariant $r$-Lagrangian of weight $\rho=0$ and of order $k=1$ with values in the
oriented Grassmannian $G^+_r$.

Let $\omega$ be an arbitrary closed $r$-form on $G_r^+$ and let the Lagrangian $L_\omega$
correspond to $\omega$ (see \eqref{diffformdensitycorrespondence}). Then the composition
$L_\omega\circ F_m$ of Lagrangians $L$ and $F_m$  is a closed covariant Lagrangian of
weight $\rho=1$ and order $k=2$. If  forms $\omega$ are not cohomologous to zero we come
to covariant Lagrangians $L_\omega\circ F$ corresponding to top-degree characteristic
classes of surfaces embedded in $\mathbb R^n$. We shall consider from this point of view
the Euler classes for $r$-dimensional surfaces embedded in $\mathbb R^m$ in two cases:
$m=n-1$, $n=2k+1$ (even-dimensional hypersurfaces in $\mathbb R^m$) and $m=2$
(two-dimensional surface, embedded in an arbitrary $\mathbb R^m$).

Case 1 ($r=m-1$). The oriented Grassmannian $G^+_{m-1}$ is simply the sphere $S^{m-1}$.
Let $\omega$ be a volume form on $S^{m-1}$:
\begin{equation}\label{gaussmap}
  \omega=\iota_E\left(r^{-m}\volx\right)=\sum
  {(-1)^{a-1}x^a dx^1\dots \widehat{dx^a}\dots dx^m\over r^{m}},
\end{equation}
where $\volx$ is the standard coordinate volume form on $\mathbb R^m$, $E=X^a{\p/\p x^a}$
is the Euler field and $r^2=(x^1)^2+\ldots+(x^m)^2$. Clearly, $F_{m-1}$ supplies the
Gauss spherical map for each oriented hypersurface $C^{m-1}\subset \R{m}$. It is not
difficult to see that the value of the covariant Lagrangian $L=L_\omega\circ F_{m-1}$
gives the Gauss curvature density. In the case $m-1=2k$, $\int L$ is the Euler class.

Case 2 ($r=2$). Let $L$ be a covariant $2$-Lagrangian in the Euclidean space $\mathbb
R^m$ such that for every two-dimensional oriented surface $C\subset \mathbb R^m$,
\begin{equation}\label{eulerclass}
   \int_C L\,\volt =\int L(x(t,\der{x}{t},\ldots)\,\volt=\int_C R\sqrt{\det g}\,\volt,
\end{equation}
where $x^a(t^i)$ is an arbitrary parametrization of the surface $C$, ($i=1,2$),
$g=(g_{ij})$ is the Riemannian metric induced on the surface $C$:
\begin{equation}
  g_{ij}=\sum {\p x^a\over\p t^i}{\p x^a\over\p t^j}
\end{equation}
and $R$ is the scalar curvature of this metric. Straightforward calculations show that
\begin{equation}
               \label{explicitgaussbonet}
                  L=\sum_{a,b}
                  {
                \left(
     x_{11}^a x_{22}^b-x_{12}^a x_{12}^b
                 \right)
                   P^{ab}
                  \over
                  \sqrt {\det g}
                  }\,,
\end{equation}
where
\begin{equation}
           \label{projector}
    P^{ab}=\delta^{ab}-x^a_i g^{ij}x^b_j
\end{equation}
is the projector on the plane orthogonal to the surface $C$, $g^{ij}$ is the tensor
inverse to the metric tensor $g_{ij}$.

We shall represent the covariant $2$-Lagrangian~\eqref{explicitgaussbonet} as the
composition of the Lagrangians $L_\omega$ and $F_2$, where $\omega$ is a closed $2$-form
on Grassmannian $G^+_2$. Consider on the Grassmannian $G_2^+$ the homogeneous coordinates
$(u^a_i)$ ($a=1,\ldots,m$, $i=1,2$) such that $u_i$ and ${u'_i}$ define the same point
(oriented 2-subspace of $\mathbb R^m$ spanned by the vectors $u_1$, $u_2$) if $
{u'_i}=t_i^j u_j$, where $t^i_j$ is a $2\times 2$ matrix with the positive determinant.
In these coordinates the covariant Lagrangian $F_2$ with values in the Grassmannian
$G_2^+$ is $F_2(x^a_1,x^a_2)=(x^a_1,x^a_2)$. Thus the covariant $2$-Lagrangian $L$
in~\eqref{explicitgaussbonet} is equal to $L_\omega\circ F_2$, where the $2$-form
$\omega$ on $G_2^+$ is given by the following formula:
\begin{equation}
         \label{formongrassmannian}
                     \omega=\sum_{a,b}
                           {
                    du^a_1 du^b_2\,
                      P^{ab}(u_1,u_2)
                         \over
                     \sqrt { \det g(u_1,u_2)}
                            },
\end{equation}
where $g(u_1,u_2)$ is the Gram matrix for the $2$-frame $u_1$, $u_2$ and
$P^{ab}=\delta^{ab}-u^a_i g^{ij} u^b_j$
is the projector on the plane orthogonal to the plane spanned by $u_1,u_2$.

The cohomology class of the closed $2$-form~\eqref{formongrassmannian} is a generator of
the group $H^2(G_2^+(\R{m}))=\mathbb R$. To get a better understanding of this form
consider the orthonormal homogeneous coordinates on $G_2^+(\R{m})$, i.e., the coordinates
$n^a_i$ of the vectors $n_1,n_2$ of an orthonormal basis of the corresponding subspace
(i.e., a point of $G^+_2$). In these coordinates the form~\eqref{formongrassmannian} has
the appearance
$$ \omega=\sum dn^a_1 dn^a_2.
$$

It seems plausible that constructions using composition of Lagrangians can be helpful for
the study of topological invariants of surfaces. This approach can be naturally
generalized to the supercase, where standard geometrical considerations are unavailable.

\def\cprime{$'$}

\end{document}